\documentclass[10pt]{article}
\usepackage{amssymb}

\newcommand{\be}{\begin{equation}}
\newcommand{\ee}{\end{equation}}
\newcommand{\bea}{\begin{eqnarray}}
\newcommand{\eea}{\end{eqnarray}}
\newcommand{\barray}{\begin{array}}
\newcommand{\earray}{\end{array}}
\newcommand{\pa}{\partial}
\newcommand{\nn}{\nonumber}
\newcommand{\bitem}{\begin{itemize}}
\newcommand{\eitem}{\end{itemize}}
\newtheorem{teo}{Theorem}[section]
\newcommand{\bt}{\begin{teo}}
\newcommand{\et}{\end{teo}}
\newtheorem{Def}{Definition}[section]
\newcommand{\bd}{\begin{Def}}
\newcommand{\ed}{\end{Def}}
\newtheorem{lem}{Lemma}[section]
\newcommand{\bl}{\begin{lem}}
\newcommand{\el}{\end{lem}}
\newtheorem{prop}{Proposition}[section]
\newcommand{\bp}{\begin{prop}}
\newcommand{\ep}{\end{prop}}
\newtheorem{cor}{Corollary}[section]
\newcommand{\bc}{\begin{cor}}
\newcommand{\ec}{\end{cor}}
\newtheorem{ex}{Example}[section]
\newcommand{\bex}{\begin{ex}}
\newcommand{\eex}{\end{ex}}
\newtheorem{rem}{Remark}[section]
\newcommand{\br}{\begin{rem}}
\newcommand{\er}{\end{rem}}

\catcode `\@=11
\@addtoreset{equation}{section}

\begin{document}

\begin{center}
{\Large \textbf{Frobenius Manifolds as a Special Class \\ of
Submanifolds in Pseudo-Euclidean
Spaces\footnote{This research was supported by the
Max-Planck-Institut f\"ur Mathematik (Bonn, Germany),
the Russian Foundation for Basic Research
(Grant No. 05-01-00170) and the Program for Supporting
Leading Scientific Schools (Grant No. NSh-4182.2006.1).}}}
\end{center}

\medskip

\begin{center}
{\large {O. I. Mokhov}}
\end{center}
\begin{center}
{Center for Nonlinear Studies,
L.D.Landau Institute for Theoretical Physics,\\
Russian Academy of Sciences,
Kosygina 2, Moscow, GSP-1, 117940, Russia\\
Department of Geometry and Topology,
Faculty of Mechanics and Mathematics, \\
M.V.Lomonosov Moscow State University,
Moscow, GSP-1, 119991, Russia\\
E-mail: mokhov@mi.ras.ru; mokhov@landau.ac.ru; mokhov@bk.ru}
\end{center}

\begin{abstract}
We introduce a very natural class of {\it potential\/} submanifolds in
pseudo-Euclidean spaces (each $N$-dimensional potential submanifold is a
special flat torsionless submanifold in a
$2N$-dimensional pseudo-Euclidean space) and prove that
each $N$-dimensional Frobenius manifold can be locally represented as an
$N$-dimensional potential submanifold. We show that
all potential submanifolds
bear natural special structures of Frobenius algebras
on their tangent spaces. These special Frobenius structures are
generated by the corresponding flat first fundamental form
and the set of the second fundamental forms of the submanifolds
(in fact, the structural constants are given
by the set of the Weingarten operators of the submanifolds).
We prove that the associativity
equations of two-dimensional topological quantum field theories
are very natural reductions of the fundamental nonlinear equations
of the theory of submanifolds in
pseudo-Euclidean spaces and define locally the class of potential
submanifolds. The problem of explicit realization
of an arbitrary concrete
Frobenius manifold as a potential submanifold
in a pseudo-Euclidean
space is reduced to solving a linear system of
second-order partial
differential equations.
For concrete Frobenius manifolds,
this realization problem can be solved explicitly
in elementary and special functions.
Moreover, we consider a nonlinear system, which is a
natural generalization of the associativity equations,
namely, the system describing all flat torsionless submanifolds
in pseudo-Euclidean spaces, and prove that this system is integrable by
the inverse scattering method. We prove that each
flat torsionless submanifold in a pseudo-Euclidean space gives
a nonlocal Hamiltonian operator of hydrodynamic type with
flat metric, a special pencil of compatible Poisson structures,
a recursion operator, infinite sets of integrals of hydrodynamic type
in involution and a natural class of integrable hierarchies,
which are all directly associated with this
flat torsionless submanifold.
In particular, using our construction of
the reduction to the associativity equations,
we obtain that each Frobenius manifold
(in point of fact, each solution of the associativity equations)
gives a natural nonlocal Hamiltonian operator of hydrodynamic
type with flat metric, a natural pencil of compatible
Poisson structures (local and nonlocal), a natural
recursion operator, natural infinite sets of
integrals of hydrodynamic type
in involution and a natural class of integrable hierarchies,
which are all directly associated with this Frobenius manifold.

\end{abstract}

{\bf Key words and phrases}: Frobenius manifold, Frobenius algebra,
symmetric algebra, $N$-parameter deformation of algebra,
submanifold in pseudo-Euclidean space,
flat submanifold, submanifold with flat normal bundle,
flat submanifold with zero torsion, indefinite metric,
associativity equations in two-dimensional topological
quantum field theories, WDVV equations, topological
quantum field theory, integrable nonlinear system, integrable hierarchy,
bi-Hamiltonian system, nonlocal Hamiltonian operator of
hydrodynamic type, system of hydrodynamic type,
compatible Poisson brackets, Poisson pencil,
system of integrals in involution, conservation law, recursion operator,
pseudo-Riemannian geometry, potential submanifold.


\section{Introduction} \label{introduction}

We prove that the associativity
equations of two-dimensional topological quantum field theories
(the Witten--Dijkgraaf--Verlinde--Verlinde equations,
see \cite{1d}--\cite{4d}) for
a function (a {\it potential\/})
$\Phi = \Phi (u^1, \ldots, u^N)$,
\be
\sum_{k = 1}^N \sum_{l = 1}^N
{\pa^3 \Phi \over \pa u^i \pa u^j \pa u^k} \eta^{kl}
{\pa^3 \Phi \over \pa u^l \pa u^m \pa u^n} =
\sum_{k = 1}^N \sum_{l = 1}^N
{\pa^3 \Phi \over \pa u^i \pa u^m \pa u^k} \eta^{kl}
{\pa^3 \Phi \over \pa u^l \pa u^j \pa u^n},  \label{ass1}
\ee
where
$\eta^{ij}$ is an arbitrary constant nondegenerate symmetric matrix,
$\eta^{ij} = \eta^{ji},$ $\eta^{ij} = {\rm const},$
$\det (\eta^{ij})\neq0$,
are very natural reductions of the fundamental nonlinear equations
of the theory of submanifolds in
pseudo-Euclidean spaces (namely, the Gauss equations,
the Codazzi equations and the Ricci equations)
and give locally a very natural class of {\it potential\/}
submanifolds in pseudo-Euclidean spaces.
Each $N$-dimensional potential submanifold is a
special flat torsionless submanifold in a
$2N$-dimensional pseudo-Euclidean space. All potential submanifolds
in pseudo-Euclidean spaces bear natural special structures of
Frobenius algebras
on their tangent spaces. These special Frobenius structures are
generated by the corresponding flat first fundamental form
and the set of the second fundamental forms of the submanifolds
(in fact, the structural constants are given
by the set of the Weingarten operators of the submanifolds).

We recall that the associativity equations (\ref{ass1}) are
consistent and integrable by the inverse scattering method,
they possess a rich set of nontrivial solutions,
and each solution
$\Phi (u^1, \ldots, u^N)$ of the
associativity equations (\ref{ass1}) gives $N$-parameter
deformations of special
Frobenius algebras (some special
commutative associative algebras equipped with
nondegenerate invariant symmetric bilinear forms) (see \cite{1d}).
Indeed, consider algebras $A (u)$ in an $N$-dimensional vector space with
the basis $e_1, \ldots, e_N$ and the multiplication (see \cite{1d})
\be
e_i \circ e_j = c^k_{ij} (u) e_k, \ \ \ \
c^k_{ij} (u) = \eta^{ks} {\pa^3 \Phi \over \pa u^s \pa u^i \pa u^j}.
\label{al1}
\ee
For all values of the parameters $u = (u^1, \ldots, u^N)$ the algebras
$A (u)$ are commutative,
$e_i \circ e_j = e_j \circ e_i,$ and the associativity condition
\be
(e_i \circ e_j) \circ e_k = e_i \circ (e_j \circ e_k)  \label{al2}
\ee
in the algebras $A (u)$ is equivalent to equations (\ref{ass1}).
The matrix $\eta_{ij}$ inverse to the matrix $\eta^{ij}$,
$\eta^{is} \eta_{sj} = \delta^i_j$, defines a
nondegenerate invariant symmetric bilinear form on the algebras
$A (u)$,
\be
\langle e_i, e_j \rangle  = \eta_{ij}, \ \ \ \
\langle e_i \circ e_j, e_k \rangle =
\langle e_i, e_j \circ e_k \rangle.  \label{al3}
\ee
Recall that locally the tangent space at every point of
any Frobenius manifold (see \cite{1d}) possesses the structure of
Frobenius algebra
(\ref{al1})--(\ref{al3}), which is determined by a solution of the
associativity equations (\ref{ass1}) and smoothly depends on the point.
We prove that
each $N$-dimensional Frobenius manifold can be locally represented as an
$N$-dimensional potential flat torsionless submanifold
in a $2N$-dimensional pseudo-Euclidean
space. The problem of explicit realization
of an arbitrary concrete $N$-dimensional
Frobenius manifold as an
$N$-dimensional potential flat torsionless submanifold
in a $2N$-dimensional pseudo-Euclidean
space is reduced to solving a linear system of
second-order partial
differential equations.
For concrete Frobenius manifolds,
this realization problem can be solved explicitly
in elementary and special functions. We shall give many
explicit important examples of these realizations in a separate paper.

Moreover, we consider a nonlinear system, which is a
natural generalization of the associativity equations (\ref{ass1}),
namely, the system describing all flat torsionless submanifolds
in pseudo-Euclidean spaces, and prove that this system is integrable by
the inverse scattering method. We also prove that each
flat torsionless submanifold in a pseudo-Euclidean space
gives a natural nonlocal Hamiltonian operator of hydrodynamic type with
flat metric, a natural special pencil of compatible Poisson structures,
a natural recursion operator, natural infinite sets of integrals of
hydrodynamic type in involution and a natural class of
integrable hierarchies, which are all directly associated with this
flat torsionless submanifold. In particular,
using our construction of the reduction to the associativity equations,
we obtain that each Frobenius manifold
(in point of fact, each solution of the associativity equations
(\ref{ass1})) gives a natural nonlocal Hamiltonian operator of hydrodynamic
type with flat metric, a natural special pencil of compatible
Poisson structures (local and nonlocal), a natural
recursion operator, natural infinite sets of
integrals of hydrodynamic type
in involution and a natural class of integrable hierarchies,
which are all directly associated with this Frobenius manifold.

\section{Frobenius algebras, Frobenius manifolds
and associativity \\equations} \label{section1a}

\subsection{Frobenius and symmetric algebras}

Recall the notion of Frobenius algebra over a field
${\mathbb{K}}$ (in this paper we consider Frobenius algebras only over
${\mathbb{R}}$ or ${\mathbb{C}}$). First of all, we must note that
there are various conventional definitions of Frobenius algebras.
In particular, sometimes in mathematical literature
a finite dimensional algebra
${\mathcal{A}}$ (with multiplication $\circ$) over a field
${\mathbb{K}}$ is called {\it Frobenius\/} if it is equipped with
a linear functional
\be
\theta : {\mathcal{A}} \rightarrow {\mathbb{K}}  \label{lf}
\ee
such that if $\theta (a \circ b) = 0$ for all
$a \in {\mathcal{A}}$, then $b = 0$.
In this case, ${\rm Ker\,} \theta$ contains no
nontrivial ideals.
It is also obvious that the bilinear form
$f (a, b) = \theta (a \circ b)$ is nondegenerate for every such linear
functional in any finite dimensional algebra.
If algebra is associative, then
we have
\be
f (a \circ b, c) = \theta ((a \circ b) \circ c) =
\theta (a \circ (b \circ c)) = f (a, b \circ c)
\ee
for all $a, b, c \in {\mathcal{A}}$ ({\it invariance\/} or
{\it associativity\/} of bilinear form).

\bd
{\rm
A bilinear form
$f : {\mathcal{A}} \times {\mathcal{A}} \rightarrow {\mathbb{K}}$
in an algebra ${\mathcal{A}}$ is called {\it invariant\/}
(or {\it associative\/}) if
\be
f (a \circ b, c) = f (a, b \circ c)
\ee
for all
$a, b, c \in {\mathcal{A}}$.
}
\ed

Consider the following conventional general definition of
Frobenius algebra.

\bd
{\rm
A finite dimensional algebra ${\mathcal{A}}$ over a field
${\mathbb{K}}$ is called {\it Frobenius\/} if it is equipped with
a nondegenerate invariant bilinear form.
}
\ed

Generally speaking, even associativity of algebra is not assumed
here (we note that often some of the following additional conditions are
included in definition of Frobenius algebras:
symmetry of invariant bilinear form,
presence of a unit in algebra, associativity of
algebra, and commutativity of algebra).

Consider an arbitrary Frobenius algebra $({\mathcal{A}}, f)$, an arbitrary
element $w \in {\mathcal{A}}$ and the corresponding linear functional
$\theta_w (a) = f (a, w)$ in ${\mathcal{A}}$. Then we have
$\theta_w (a \circ b) = f (a \circ b, w) =
f (a, b \circ w)$. Therefore, if $\theta_w (a \circ b) =
f (a, b \circ w) = 0$ for all
$a \in {\mathcal{A}}$, then $b \circ w = 0$. If $w$ is an element
of algebra ${\mathcal{A}}$ such that $b \circ w = 0$
implies $b = 0$,
then $\theta_w (a)$ is a linear functional of type
(\ref{lf}) and ${\rm Ker\,} \theta_w$ contains no ideals.
For example, if algebra contains a unit $e$, then
the unit $e$ gives a linear functional of type
(\ref{lf}), $\theta_e (a) = f (a, e)$, and
${\rm Ker\,} \theta_e$ contains no ideals. Moreover,
for any algebra with a unit $e$, any invariant bilinear form
$f$ is completely generated by the linear functional
$\theta_e (a) = f (a, e)$, since $f (a, b) = f (a, b \circ e) =
f (a \circ b, e) = \theta_e (a \circ b)$.

\bex {\bf Matrix algebra $M_n ({\mathbb{K}})$}.

{\rm
Consider the algebra $M_n ({\mathbb{K}})$ of $n \times n$ matrices over
a field ${\mathbb{K}}\,$, the linear functional (trace of matrices)
$$\theta (a) = {\rm Tr\,} (a), \ a \in M_n ({\mathbb{K}}),$$
and the bilinear form $f (a, b) = \theta (a b)$.
The bilinear form is invariant,
since the matrix algebra is associative. It is easy to prove that
the bilinear form is nondegenerate, and
$(M_n ({\mathbb{K}}), f)$ is a noncommutative associative Frobenius
algebra with a unit over ${\mathbb{K}}\,$.
Note that the bilinear form
$f (a, b) = \theta (a b)$ is symmetric, $\theta (a b) =
\theta (b a)$. Recall that a finite dimensional associative algebra
with a unit over a field
${\mathbb{K}}$ is called {\it symmetric\/} if it is equipped with
a symmetric nondegenerate associative bilinear form (see
\cite{curtis}).
Therefore, $(M_n ({\mathbb{K}}), f)$ is a symmetric algebra.
}
\eex

\bex  {\bf Group algebra ${\mathbb{K}}G$}.

{\rm
Let $G$ be a finite group. Consider the group algebra ${\mathbb{K}}G$
over a field ${\mathbb{K}}\,$,
$${\mathbb{K}}G =
\{a \ | \ a = \sum_{g \in G} \alpha_g g, \  \alpha_g \in {\mathbb{K}} \}.$$
${\mathbb{K}}G$ is an associative algebra with a unit
over ${\mathbb{K}}\,$.
Let $e$ be the unit of the group $G$. Consider the linear
functional
$$ \theta (a) = \alpha_e (a),
\ \ a = \sum_{g \in G} \alpha_g (a) g \in {\mathbb{K}}G, \
\ \alpha_g (a) \in {\mathbb{K}},$$
and the bilinear form
$f (a, b) = \theta (a b)$. The bilinear form is invariant,
since the group algebra is associative. It is easy to prove that
the bilinear form is
nondegenerate. Indeed, we have
$$f (g^{-1}, a) = \theta (g^{-1} a) = \alpha_g (a)$$
for all $g \in G$. Therefore, if $f (g, a) = \theta (g a) = 0$ for all
$g \in G$, then $\alpha_g (a) = 0$ for all $g \in G$, i.e., $a = 0$.
Hence the bilinear form $f$ is nondegenerate, and
$({\mathbb{K}}G, f)$ is a noncommutative associative Frobenius
algebra with a unit over ${\mathbb{K}}$ (it is commutative only
for Abelian groups).
Note that the bilinear form
$f (a, b) = \theta (a b)$ is symmetric for any group $G$, $\theta (a b) =
\theta (b a)$. Therefore, $({\mathbb{K}}G, f)$ is a symmetric algebra.
}
\eex

\subsection{Frobenius manifolds}

Consider an $N$-dimensional pseudo-Riemannian manifold $M$
with a metric $g$ and a structure of Frobenius algebra
$(T_u M, \circ, g)$,
$T_u M \times T_u M \stackrel{\circ}{\rightarrow} T_u M$,
on each tangent space $T_u M$ at any point $u \in M$ smoothly
depending on the point such that the metric $g$ is the
corresponding nondegenerate
invariant symmetric bilinear form on each tangent space $T_u M$,
$g (X \circ Y, Z) = g (X, Y \circ Z),$
where $X, Y$ and $ Z$ are arbitrary vector fields on $M$.

This class of pseudo-Riemannian manifolds
equipped with Frobenius structures could be naturally
called Frobenius, but in this paper we shall consider well-known and
generally accepted Dubrovin's
definition of Frobenius manifolds \cite{1d}, which is motivated by
two-dimensional topological quantum field theories and quantum cohomology
and imposes very severe additional constraints on Frobenius structures
of Frobenius manifolds.

\bd {\rm (Dubrovin \cite{1d})

An $N$-dimensional pseudo-Riemannian manifold $M$
with a metric $g$ and a structure of Frobenius algebra
$(T_u M, \circ, g)$,
$T_u M \times T_u M \stackrel{\circ}{\rightarrow} T_u M$,
on each tangent space $T_u M$ at any point $u \in M$ smoothly
depending on the point
is called {\it Frobenius\/} if

(1) the metric $g$ is a nondegenerate
invariant symmetric bilinear form on each tangent space $T_u M$,
\be
g (X \circ Y, Z) = g (X, Y \circ Z), \label{f1a}
\ee

(2) the Frobenius algebra is commutative,
\be
X \circ Y = Y \circ X
\ee
for all vector fields $X$ and $Y$ on $M$,

(3) the Frobenius algebra is associative,
\be
(X \circ Y) \circ Z = X \circ (Y \circ Z)
\ee
for all vector fields $X, Y$ and $Z$ on $M$,

(4) the metric $g$ is flat,

(5) $A (X, Y, Z) = g (X \circ Y, Z)$ is a symmetric tensor on $M$
(it is obvious that, by virtue of (1) and (2), we have
$g (X \circ Y, Z) = g (X, Y \circ Z) = g (Y \circ Z, X) =
g (Y, Z \circ X) = g (Z \circ X, Y) = g (Z, X \circ Y) =
g (Z, Y \circ X) = g (Z \circ Y, X) =
g (X, Z \circ Y) = g (X \circ Z, Y) = g (Y, X \circ Z) =
g (Y \circ X, Z)$) such
that the tensor $(\nabla_W A) (X, Y, Z)$ is symmetric with respect
to all vector fields $X, Y, Z$ and $W$ on  $M$ ($\nabla$ is the covariant
differentiation generated by the
Levi-Civita connection of the metric $g$),

(6) the Frobenius algebra possesses a unit, and
the unit vector field $U$,
for which $X \circ U = U \circ X = X$ for each
vector field $X$ on $M$, is covariantly constant, i.e.,
\be
\nabla U = 0,
\ee
where $\nabla$ is the covariant
differentiation generated by the
Levi-Civita connection of the metric $g$,

(7) the manifold $M$ is equipped with
a vector field $E$ ({\it Euler vector field}) such that
\be
\nabla \nabla E = 0,
\ee
\be
{\mathcal {L}}_E (X \circ Y) - ({\mathcal {L}}_E X ) \circ Y -
X \circ ({\mathcal {L}}_E Y) = X \circ Y,
\ee
\be
{\mathcal {L}}_E \, g (X, Y)
- g ({\mathcal {L}}_E X,  Y) -
g (X, {\mathcal {L}}_E Y) = K\, g (X, Y),
\ee
\be
{\mathcal {L}}_E U = - U,
\ee
where $K$ is an arbitrary fixed constant,
${\mathcal {L}}_E$ is the Lie derivative along the Euler vector field,
and $\nabla$ is the covariant
differentiation generated by the
Levi-Civita connection of the metric $g$.
}
\ed

A beautiful theory of these very special
Frobenius structures and Frobenius manifolds and many
important examples
were constructed by Dubrovin in connection with
two-dimensional topological quantum field theories
and quantum cohomology \cite{1d}.
No doubt that these very special Frobenius structures
and Frobenius manifolds should be called
Dubrovin's. A lot of very important examples of Frobenius manifolds
arises in the theory of Gromov--Witten invariants, the quantum cohomology,
the singularity theory, the enumerative geometry, the topological
field theories and the
modern differential geometry, mathematical and theoretical physics.

In this paper we describe a very natural special class of
submanifolds in pseudo-Euclidean spaces bearing natural
Frobenius structures satisfying
the conditions (1)--(5), namely, the class of potential
submanifolds. Moreover, we show that
each manifold satisfying
the conditions (1)--(5) can be locally realized as a potential submanifold
in a pseudo-Euclidean space
\cite{2m}--\cite{4m}. For any concrete
Frobenius structure satisfying
the conditions (1)--(5) and for any given Frobenius
manifold, the corresponding realization
problem is reduced to solving a system of
linear second-order partial differential equations.

\subsection{Associativity equations}

Consider an arbitrary manifold satisfying
the conditions (1)--(5). Let $u = (u^1, \ldots, u^N)$
be arbitrary flat coordinates of the flat metric $g$.
In flat local coordinates,
the metric $g (u)$ is a constant nondegenerate
symmetric matrix $\eta_{ij}$, $\eta_{ij} = \eta_{ji},$
$\det (\eta_{ij}) \neq 0,$ $\eta_{ij} = \  {\rm const}$,
$g (X, Y) = \eta_{ij} X^i (u) Y^j (u)$.

In these flat local coordinates,
for structural functions $c^i_{jk} (u)$
of the Frobenius structure on the manifold,
$$X \circ Y = W, \ \ \ W^i (u) = c^i_{jk} (u) X^j (u) Y^k (u),$$
and for the symmetric tensor $A_{ijk} (u)$, we have
\bea
&&
A (X, Y, Z) = A_{ijk} (u) X^i (u) Y^j (u) Z^k (u)
= g (X \circ Y, Z) = \nn\\
&&
= g (W, Z) = \eta_{ij} W^i (u) Z^j (u)
= \eta_{ij} c^i_{kl} (u) X^k (u) Y^l (u) Z^j (u).\nn
\eea
Therefore,
\be
A_{ijk} (u) = \eta_{sk} c^s_{ij} (u). \label{str}
\ee

According to (5) $(\nabla_l A_{ijk}) (u)$ is a symmetric tensor,
i.e., in the flat local coordinates we also have
$${\pa A_{ijk} \over \pa u^l} = {\pa A_{ijl} \over \pa u^k}.$$
Hence there locally exist functions $B_{ij} (u)$ such that
$$A_{ijk} (u) = {\pa B_{ij} \over \pa u^k}.$$
We can consider that the matrix $B_{ij} (u)$ is symmetric,
$B_{ij} (u) = B_{ji} (u)$. Indeed, if
$$A_{ijk} (u) = {\pa \widetilde{B}_{ij} \over \pa u^k},$$
then
$${\pa \widetilde{B}_{ij} \over \pa u^k} =
{\pa \widetilde{B}_{ji} \over \pa u^k}$$
for any $k$, since the tensor $A_{ijk} (u)$ is symmetric.
Hence, $\widetilde{B}_{ij} (u) =
\widetilde{B}_{ji} (u) + C_{ij}$, where
$C_{ij} = \ {\rm const},$ $C_{ij} = - C_{ji}.$
Thus, if we take $B_{ij} (u) = \widetilde{B}_{ij} (u) - (1 / 2) C_{ij},$
then
$B_{ij} (u) = B_{ji} (u)$ and
$$A_{ijk} (u) = {\pa B_{ij} \over \pa u^k}.$$
Since the tensor $A_{ijk} (u)$ is symmetric, we have also
$${\pa B_{ij} \over \pa u^k} = {\pa B_{ik} \over \pa u^j}.$$
Hence there locally exist functions $F_i (u)$ such that
$$B_{ij} (u) = {\pa F_i \over \pa u^j}.$$
Since the matrix $B_{ij} (u)$ is symmetric, we have
$${\pa F_i \over \pa u^j} = {\pa F_j \over \pa u^i}.$$
Hence there locally exist a function
(a {\it potential\/}) $\Phi (u)$ such that
$$F_i (u) = {\pa \Phi \over \pa u^i}.$$
Thus,
$$A_{ijk} (u) = {\pa B_{ij} \over \pa u^k} =
{\pa^2 F_i \over \pa u^j \pa u^k} =
{\pa^3 \Phi \over \pa u^i \pa u^j \pa u^k}.$$

From (\ref{str}) for the structural functions $c^i_{jk} (u)$ we have
\be
c^i_{jk} (u) = \eta^{is} A_{sjk} (u) =
\eta^{is} {\pa^3 \Phi \over \pa u^s \pa u^j \pa u^k}, \label{str2}
\ee
where the matrix $\eta^{ij}$ is inverse to the matrix $\eta_{ij}$,
$\eta^{is} \eta_{sj} = \delta^i_j$.

For any values of the parameters $u = (u^1, \ldots, u^N)$,
the structural functions (\ref{str2}) give a commutative
Frobenius algebra
\be
\pa_i \circ \pa_j = c^k_{ij} (u) \pa_k =
\eta^{ks} {\pa^3 \Phi \over \pa u^s \pa u^i \pa u^j} \pa_k \label{frob}
\ee
equipped with a symmetric invariant
nondegenerate bilinear form
\be
\langle \pa_i, \pa_j \rangle = \eta_{ij} \label{form}
\ee
for any constant nondegenerate
symmetric matrix $\eta_{ij}$ and for any function $\Phi (u)$,
but, generally speaking, this algebra is not associative.
All the conditions (1)--(5) except the associativity condition (3) are
obviously satisfied for all
these $N$-parameter deformations of nonassociative Frobenius
algebras.

The associativity condition (3) is equivalent to
a nontrivial overdetermined system of nonlinear partial
differential equations for the potential $\Phi (u)$,
\be
\sum_{k = 1}^N \sum_{l = 1}^N
{\pa^3 \Phi \over \pa u^i \pa u^j \pa u^k} \eta^{kl}
{\pa^3 \Phi \over \pa u^l \pa u^m \pa u^n} =
\sum_{k = 1}^N \sum_{l = 1}^N
{\pa^3 \Phi \over \pa u^i \pa u^m \pa u^k} \eta^{kl}
{\pa^3 \Phi \over \pa u^l \pa u^j \pa u^n},  \label{ass1a}
\ee
which is well known as the associativity
equations of two-dimensional topological quantum field theories
(the Witten--Dijkgraaf--Verlinde--Verlinde or the WDVV equations,
see \cite{1d}--\cite{4d}); it is
consistent, integrable by the inverse scattering method
and possesses a rich set of nontrivial solutions (see \cite{1d}).

It is obvious that each solution
$\Phi (u^1, \ldots, u^N)$ of the
associativity equations (\ref{ass1a}) gives $N$-parameter deformations of
commutative associative Frobenius algebras
(\ref{frob})
equipped with
nondegenerate invariant symmetric bilinear forms (\ref{form}). These
Frobenius structures satisfy to all the conditions (1)--(5).

Further in this paper we show that the associativity
equations (\ref{ass1a})
are very natural reductions of the fundamental nonlinear equations
of the theory of submanifolds in
pseudo-Euclidean spaces and give a natural class of {\it potential\/}
flat torsionless submanifolds \cite{2m}--\cite{4m}.
All potential flat torsionless submanifolds
in pseudo-Euclidean spaces bear natural structures of Frobenius algebras
(\ref{frob}), (\ref{form}) on their tangent spaces.
These Frobenius structures are
generated by the corresponding flat first fundamental form
and the set of the second fundamental forms of the submanifolds.

\section{Gauss, Codazzi, and Ricci equations and
Bonnet theorem in \\the theory of submanifolds in
Euclidean spaces} \label{section1}

\subsection{Submanifolds in Euclidean spaces}

Let us consider an arbitrary smooth $N$-dimensional submanifold
$M^N$ in an $(N+L)$-di\-mensional Euclidean space $E^{N + L}$,
$M^N \subset E^{N + L}$,
and introduce the standard classical notation. Let the submanifold
$M^N$ be given locally by a smooth vector function
$ r (u^1, \ldots, u^N)$ of $N$ independent variables
$(u^1, \ldots, u^N)$ (some independent parameters on the submanifold),
$ r (u^1, \ldots, u^N) =
(z^1 (u^1, \ldots, u^N), \ldots, z^{N + L} (u^1, \ldots, u^N)),$
where $(z^1, \ldots, z^{N + L})$ are Cartesian coordinates in the
Euclidean space $E^{N + L}$, $(z^1, \ldots, z^{N + L}) \in E^{N + L}$,
$(u^1, \ldots, u^N)$ are local coordinates (parameters) on $M^N$,
${\rm rank\,} (\pa z^i/\pa u^j) = N$ (here $1 \leq i \leq N + L,$
$1 \leq j \leq N$).
Then ${\pa r/ \pa u^i} = r_{u^i},$ $1 \leq i \leq N,$ are
tangent vectors at any point $u = (u^1, \ldots, u^N)$ on $M^N$.
Let ${\bf N}_u$ be the
normal space of the submanifold
$M^N$ at an arbitrary point
$u = (u^1, \ldots, u^N)$ on $M^N$,
${\bf N}_u = \langle n_1, \ldots, n_L \rangle$, where
$n_{\alpha}$, $1 \leq \alpha \leq L,$
is an orthonormalized basis of the normal space
(orthonormalized normals),
$(n_{\alpha}, r_{u^i}) = 0,$ $1 \leq \alpha \leq L,$
$1 \leq i \leq N,$
$(n_{\alpha}, n_{\beta}) = \delta_{\alpha \beta},$
$1 \leq \alpha, \beta \leq L$.
Then ${\bf I} = d s^2 = g_{ij} (u) d u^i d u^j,$
$g_{ij} (u) = (r_{u^i}, r_{u^j}),$
is the first fundamental form, and
${\bf II}_{\alpha} = \omega_{\alpha, ij} (u) d u^i d u^j,$
$\omega_{\alpha, ij} (u) = (n_{\alpha}, r_{u^i u^j}),$
$1 \leq \alpha \leq L,$ are the second fundamental forms
of the submanifold $M^N$.

\subsection{Torsion forms of submanifolds in Euclidean spaces}

Since
the set of vectors
$(r_{u^1} (u), \ldots, r_{u^N} (u), n_1 (u), \ldots, n_L (u))$
forms a basis in $E^{N + L}$ at each point of the submanifold $M^N$, we
can decompose each of the vectors $n_{\alpha, u^i} (u),$
$1 \leq \alpha \leq L,$
$1 \leq i \leq N,$ with respect to this basis, namely,
$$n_{\alpha, u^i} (u) =
\sum_{k = 1}^N  A^k_{\alpha, i} (u) r_{u^k} (u) +
\sum_{\beta = 1}^L  \varkappa_{\alpha \beta, i} (u) n_{\beta} (u),$$
where $A^k_{\alpha, i} (u)$
and $\varkappa_{\alpha \beta, i} (u)$ are some coefficients
depending on $u$ (the {\it Weingarten decomposition\/}).
It is easy to prove that
$A^k_{\alpha, i} (u) = - \omega_{\alpha, ij} (u) g^{jk} (u),$
where $g^{jk} (u)$ is the contravariant metric inverse to the
first fundamental form $g_{ij} (u)$, $g^{is} (u) g_{sj} (u) = \delta^i_j.$
The coefficients $\varkappa_{\alpha \beta, i} (u)$ are called
the {\it torsion coefficients of the submanifold $M^N$},
$\varkappa_{\alpha \beta, i} (u) = (n_{\alpha, u^i} (u), n_{\beta} (u)).$
It is also easy to prove that the coefficients
$\varkappa_{\alpha \beta, i} (u)$ are skew-symmetric with respect to
the indices $\alpha$ and $\beta$,
$\varkappa_{\alpha \beta, i} (u) = - \varkappa_{\beta \alpha, i} (u)$,
and form covariant tensors (1-forms) with respect to the index $i$ on
the submanifold $M^N$.
The 1-forms $\varkappa_{\alpha \beta, i} (u) d u^i$ are called
the {\it torsion forms of the submanifold $M^N$}.

\subsection{Fundamental nonlinear equations in the theory
of submanifolds in \\Euclidean spaces}

It is well known that for each submanifold $M^N$ the forms $g_{ij} (u)$,
$\omega_{\alpha, ij} (u)$ and $\varkappa_{\alpha \beta, i} (u)$
satisfy the Gauss equations, the Codazzi equations and the Ricci equations,
which are the fundamental equations of the theory of submanifolds.
In our case, the Gauss equations have the form
\be
R_{ijkl} (u) = \sum_{\alpha = 1}^L
\left (\omega_{\alpha, jl} (u) \omega_{\alpha, ik} (u) -
\omega_{\alpha, jk} (u) \omega_{\alpha, il} (u) \right ), \label{g1}
\ee
where $R_{ijkl} (u)$ is the Riemannian curvature tensor of
the first fundamental form $g_{ij} (u)$, the Codazzi equations
have the form
\be
\nabla_k \omega_{\alpha, ij} (u) -
\nabla_j \omega_{\alpha, ik} (u) =
\sum_{\beta = 1}^L (\varkappa_{\alpha \beta, k} (u)
\omega_{\beta, ij} (u) - \varkappa_{\alpha \beta, j} (u)
\omega_{\beta, ik} (u)), \label{c1}
\ee
where $\nabla_k$ is the covariant differentiation generated
by the Levi-Civita connection of the first fundamental form
$g_{ij} (u)$, and the Ricci equations have the form
\bea
&&
\nabla_k \varkappa_{\alpha \beta, i} (u) -
\nabla_i \varkappa_{\alpha \beta, k} (u) +
\sum_{\gamma = 1}^L \left (\varkappa_{\alpha \gamma, i} (u)
\varkappa_{\gamma \beta, k} (u) - \varkappa_{\alpha \gamma, k} (u)
\varkappa_{\gamma \beta, i} (u) \right ) + \nn\\
&&
+ \sum_{l = 1}^N \sum_{j = 1}^N  g^{lj} (u) \,
\left (\omega_{\alpha, kl} (u) \omega_{\beta, ji} (u) -
\omega_{\alpha, il} (u)
\omega_{\beta, jk} (u) \right ) = 0. \label{r1}
\eea

\subsection{Bonnet theorem for submanifolds in Euclidean spaces}

{\bf Theorem} (Bonnet).
{\it Let $K^N$ be an arbitrary smooth $N$-dimensional Riemannian
manifold with a metric
$g_{ij} (u) du^i du^j$. Let some 2-forms
$\omega_{\alpha, ij} (u) du^i du^j$, $1 \leq \alpha \leq L,$
and some 1-forms $\varkappa_{\alpha \beta, i} (u) du^i$,
$1 \leq \alpha, \beta \leq L,$ be given in a simply connected
domain of the manifold $K^N$. If $\omega_{\alpha, ij} (u) =
\omega_{\alpha, ji} (u)$, $\varkappa_{\alpha \beta, i} (u) =
- \varkappa_{\beta \alpha, i} (u)$, and the Gauss
equations {\rm (\ref{g1})},
the Codazzi equations {\rm (\ref{c1})} and the Ricci equations
{\rm (\ref{r1})} are satisfied for
the forms $g_{ij} (u)$,
$\omega_{\alpha, ij} (u)$ and $\varkappa_{\alpha \beta, i} (u)$,
then there exists a unique {\rm (}up to motions{\rm )} smooth
$N$-dimensional submanifold $M^N$
in an $(N + L)$-dimensional Euclidean space
$E^{N + L}$ with the first fundamental form
$d s^2 = g_{ij} (u) du^i du^j$, the second fundamental forms
$\omega_{\alpha, ij} (u) d u^i d u^j$ and the torsion forms
$\varkappa_{\alpha \beta, i} (u) d u^i$.}

Similar fundamental equations and the Bonnet theorem hold for all
{\it totally nonisotropic submanifolds in pseudo-Euclidean spaces\/}
(we recall that if we have a submanifold in an
arbitrary pseudo-Euclidean space $E^m_n$,
then the metric induced on the submanifold from the ambient
pseudo-Euclidean space $E^m_n$ is nondegenerate
if and only if this submanifold is totally nonisotropic, i.e.,
it is not tangent to isotropic cones of the ambient
pseudo-Euclidean space $E^m_n$ at its points).

\section{Description of flat submanifolds with zero torsion in
\\pseudo-Euclidean spaces} \label{section2}

\subsection{Submanifolds with zero torsion in pseudo-Euclidean spaces}

Let us consider totally nonisotropic smooth  $N$-dimensional
submanifolds with {\it zero torsion\/} in an
arbitrary $(N + L)$-dimensional
pseudo-Euclidean space, i.e.,
all torsion forms of submanifolds of this class vanish,
$\varkappa_{\alpha \beta, i} (u) = 0$.
In the normal spaces ${\bf N}_u$,
we also use the bases $n_{\alpha}$, $1 \leq \alpha \leq L,$
with arbitrary admissible constant Gram matrices $\mu_{\alpha \beta}$,
$(n_{\alpha}, n_{\beta}) = \mu_{\alpha \beta}$,
$\mu_{\alpha \beta}= {\rm const}$,
$\mu_{\alpha \beta} = \mu_{\beta \alpha}$, $\det (\mu_{\alpha \beta})
\neq 0$ (the signature of the metric $\mu_{\alpha \beta}$ is completely
determined by the
signature of the first fundamental form of the corresponding submanifold
and the signature of the corresponding ambient
pseudo-Euclidean space).

For torsionless $N$-dimensional submanifolds in an
arbitrary $(N + L)$-dimensional
pseudo-Euclidean space, we obtain the following
system of fundamental equations:
the Gauss equations
\be
R_{ijkl} (u) = \sum_{\alpha = 1}^L \sum_{\beta = 1}^L
\mu^{\alpha \beta}
(\omega_{\alpha, ik} (u) \omega_{\beta, jl} (u) -
\omega_{\alpha, il} (u) \omega_{\beta, jk} (u)),  \label{g2}
\ee
where $\mu^{\alpha \beta}$ is inverse to the matrix
$\mu_{\alpha \beta}$, $\mu^{\alpha \gamma} \mu_{\gamma \beta} =
\delta^{\alpha}_{\beta}$, the Codazzi equations
\be
\nabla_k \omega_{\alpha, ij} (u) =
\nabla_j \omega_{\alpha, ik} (u),   \label{c2}
\ee
and the Ricci equations
\be
\sum_{i = 1}^N \sum_{j = 1}^N
g^{ij} (u) \, (\omega_{\alpha, ik} (u) \omega_{\beta, jl} (u)
- \omega_{\alpha, il} (u)
\omega_{\beta, jk} (u)) = 0.  \label{r2}
\ee

\subsection{Second fundamental forms of
flat torsionless submanifolds
in \\pseudo-Euclidean spaces and Hessians}

Now let $g_{ij} (u)$ be a flat metric, i.e., we consider flat
torsionless $N$-dimensional submanifolds $M^N$
in an $(N + L)$-dimensional
pseudo-Euclidean space. Then we can consider that
$u = (u^1, \ldots, u^N)$ are certain flat coordinates of the
metric $g_{ij} (u)$  on $M^N$. In flat coordinates, the metric is
a constant nondegenerate symmetric matrix $\eta_{ij}$,
$\eta_{ij} = \eta_{ji},$ $\eta_{ij} = {\rm const},$
$\det (\eta_{ij}) \neq 0$, and
the Codazzi equations (\ref{c2}) have the form
\be
{\pa \omega_{\alpha, ij} \over \pa u^k} =
{\pa \omega_{\alpha, ik} \over \pa u^j}.
\ee
Therefore, there locally exist some functions $\chi_{\alpha, i} (u),$
$1 \leq \alpha \leq L,$ $1 \leq i \leq N,$ such that
\be
\omega_{\alpha, ij} (u) = {\pa \chi_{\alpha, i} \over \pa u^j}.
\ee
From symmetry of the second fundamental forms
$\omega_{\alpha, ij} (u) = \omega_{\alpha, ji} (u)$,
we have
\be
{\pa \chi_{\alpha, i} \over \pa u^j} =
{\pa \chi_{\alpha, j} \over \pa u^i}.
\ee
Therefore, there locally exist some functions $\psi_{\alpha} (u),$
$1 \leq \alpha \leq L,$ such that
\be
\chi_{\alpha, i} (u) = {\pa \psi_{\alpha} \over \pa u^i}, \ \ \ \
\omega_{\alpha, ij} (u) = {\pa^2 \psi_{\alpha} \over \pa u^i \pa u^j}.
\ee
We have thus proved the following important lemma.
\bl {\rm \cite{3m}, \cite{4m}}   \label{lem}
All the second fundamental forms of each
flat torsionless submanifold
in a pseudo-Euclidean space are Hessians in any flat coordinates
in any simply connected domain on the
submanifold.
\el

\subsection{Fundamental nonlinear equations for flat
torsionless submanifolds in \\pseudo-Euclidean spaces}

It follows from Lemma \ref{lem} that in any flat coordinates,
the Gauss equations (\ref{g2})
have the form
\be
\sum_{\alpha = 1}^L \sum_{\beta = 1}^L
\mu^{\alpha \beta}
\left ({\pa^2 \psi_{\alpha} \over \pa u^i \pa u^k}
{\pa^2 \psi_{\beta} \over \pa u^j \pa u^l} -
{\pa^2 \psi_{\alpha} \over \pa u^i \pa u^l}
{\pa^2 \psi_{\beta} \over \pa u^j \pa u^k} \right ) = 0, \label{g3}
\ee
and the Ricci equations (\ref{r2}) have the form
\be
\sum_{i = 1}^N \sum_{j = 1}^N
\eta^{ij}
\left ({\pa^2 \psi_{\alpha} \over \pa u^i \pa u^k}
{\pa^2 \psi_{\beta} \over \pa u^j \pa u^l} -
{\pa^2 \psi_{\alpha} \over \pa u^i \pa u^l}
{\pa^2 \psi_{\beta} \over \pa u^j \pa u^k} \right ) = 0,  \label{r3}
\ee
where $\eta^{ij}$ is inverse to the matrix
$\eta_{ij}$, $\eta^{is} \eta_{sj} =
\delta^i_j$.

\bt {\rm \cite{2m}--\cite{4m}}
The class of $N$-dimensional flat torsionless submanifolds
in $(N + L)$-dimensional
pseudo-Euclidean spaces is described
{\rm (}in flat coordinates{\rm )} by the
system of nonlinear equations {\rm (\ref{g3}), (\ref{r3})} for
functions $\psi_{\alpha} (u),$ $1 \leq \alpha \leq L.$
Here, $\eta^{ij}$ and $\mu^{\alpha \beta}$
are arbitrary constant nondegenerate symmetric matrices,
$\eta^{ij} = \eta^{ji},$ $\eta^{ij} = {\rm const},$
$\det (\eta^{ij}) \neq 0$,
$\mu^{\alpha \beta}= {\rm const}$,
$\mu^{\alpha \beta} = \mu^{\beta \alpha}$,
$\det (\mu^{\alpha \beta}) \neq 0${\rm ;}
the signature of the ambient $(N + L)$-dimensional
pseudo-Euclidean space is the sum of the signatures of
the metrics $\eta^{ij}$ and $\mu^{\alpha \beta}${\rm ;}
${\bf I} = d s^2 = \eta_{ij} d u^i d u^j$ is the
first fundamental form, where $\eta_{ij}$ is inverse to the matrix
$\eta^{ij}$, $\eta^{is} \eta_{sj} =
\delta^i_j$, and
${\bf II}_{\alpha} =
(\pa^2 \psi_{\alpha} / (\pa u^i \pa u^j)) d u^i d u^j,$
$1 \leq \alpha \leq L,$ are the second fundamental forms
given by the Hessians of the
functions  $\psi_{\alpha} (u),$ $1 \leq \alpha \leq L$,
for the corresponding flat torsionless submanifold.
\et

According to the Bonnet theorem,
any solution $\psi_{\alpha} (u),$ $1 \leq \alpha \leq L,$
of the nonlinear system (\ref{g3}), (\ref{r3})
determines a unique (up to motions) $N$-dimensional flat torsionless
submanifold of
the corresponding $(N + L)$-dimensional pseudo-Euclidean space with the
first fundamental form $\eta_{ij} d u^i d u^j$ and the
second fundamental forms
$\omega_{\alpha} (u) =
(\pa^2 \psi_{\alpha} / (\pa u^i \pa u^j))
d u^i d u^j$, $1 \leq \alpha \leq L$, given by the Hessians of the
functions  $\psi_{\alpha} (u),$ $1 \leq \alpha \leq L$.
It is obvious that we can always add arbitrary terms linear in
the coordinates $(u^1, \ldots, u^N)$ to any solution of the system
(\ref{g3}), (\ref{r3}), but the set of the second fundamental forms
and the corresponding submanifold
will be the same. Moreover, any two sets of
the second fundamental forms of the form $\omega_{\alpha, ij} (u) =
\pa^2 \psi_{\alpha} / (\pa u^i \pa u^j)$, $1 \leq \alpha \leq L$,
coincide if and only if the corresponding
functions $\psi_{\alpha} (u),$ $1 \leq \alpha \leq L$,
coincide up to terms linear in the coordinates; hence
we must not distinguish solutions of
the nonlinear system (\ref{g3}), (\ref{r3}) up to
terms linear in the coordinates $(u^1, \ldots, u^N)$.

\section{Integrability of the nonlinear equations for
flat torsionless \\submanifolds
in pseudo-Euclidean spaces} \label{section2a}

\subsection{Linear problem with parameters
for the nonlinear equations describing \\all
flat torsionless submanifolds in pseudo-Euclidean spaces}

Consider the following linear problem with parameters for
vector functions $\pa a(u) / \pa u^i$, $1 \leq i \leq N,$
and $b_{\alpha} (u)$, $1 \leq \alpha \leq L$:
\be
{\pa^2 a \over \pa u^i \pa u^j} = \lambda \, \mu^{\alpha \beta}
\omega_{\alpha, ij} (u)  b_{\beta} (u),\ \ \ \
{\pa b_{\alpha} \over \pa u^i} =
\rho \, \eta^{kj} \omega_{\alpha, ij} (u) {\pa a \over \pa u^k}, \label{t1}
\ee
where $\eta^{ij}$, $ 1 \leq i, j \leq N,$ and $\mu^{\alpha \beta},$
$1 \leq \alpha, \beta \leq L$,
are arbitrary constant nondegenerate symmetric matrices,
$\eta^{ij} = \eta^{ji},$ $\eta^{ij} = {\rm const},$
$\det (\eta^{ij}) \neq 0$,
$\mu^{\alpha \beta}= {\rm const}$,
$\mu^{\alpha \beta} = \mu^{\beta \alpha}$, $\det (\mu^{\alpha \beta})
\neq 0$;
$\lambda$ and $\rho$ are arbitrary constants (parameters) \cite{4m}.
Of course, only one of the parameters is essential (but it is really
essential).
It is obvious that
the coefficients
$\omega_{\alpha, ij} (u)$, $1 \leq \alpha \leq L,$ here must be
symmetric matrix functions, $\omega_{\alpha, ij} (u) =
\omega_{\alpha, ji} (u)$.

The consistency conditions for the linear system (\ref{t1})
are equivalent to the nonlinear system (\ref{g3}), (\ref{r3})
describing the class of $N$-dimensional flat torsionless submanifolds
in $(N + L)$-dimensional pseudo-Euclidean spaces. Indeed,
we have
\bea
&&
{\pa^3 a \over \pa u^i \pa u^j \pa u^k} = \lambda \, \mu^{\alpha \beta}
{\pa \omega_{\alpha, ij} \over \pa u^k}  b_{\beta} (u) +
\lambda \, \mu^{\alpha \beta}
\omega_{\alpha, ij} (u)  {\pa b_{\beta} \over \pa u^k} = \nn\\
&&
= \lambda \, \mu^{\alpha \beta}
{\pa \omega_{\alpha, ij} \over \pa u^k}  b_{\beta} (u) +
\lambda \, \mu^{\alpha \beta}
\omega_{\alpha, ij} (u)
\rho \, \eta^{ls} \omega_{\beta, ks} (u) {\pa a \over \pa u^l} = \nn\\
&&
= \lambda \, \mu^{\alpha \beta}
{\pa \omega_{\alpha, ik} \over \pa u^j}  b_{\beta} (u) +
\lambda \, \mu^{\alpha \beta}
\omega_{\alpha, ik} (u)
\rho \, \eta^{ls} \omega_{\beta, js} (u) {\pa a \over \pa u^l},
\eea
whence we obtain
\be
{\pa \omega_{\alpha, ij} (u) \over \pa u^k} =
{\pa \omega_{\alpha, ik} (u) \over \pa u^j} \label{a1}
\ee
and
\be
\mu^{\alpha \beta}
\omega_{\alpha, ij} (u)
\omega_{\beta, ks} (u) =
\mu^{\alpha \beta}
\omega_{\alpha, ik} (u)
\omega_{\beta, js} (u).  \label{b1}
\ee
Moreover,
\bea
&&
{\pa^2 b_{\alpha} \over \pa u^i \pa u^l} =
\rho \, \eta^{kj} {\pa \omega_{\alpha, ij} \over \pa u^l}
 {\pa a \over \pa u^k} +
\rho \, \eta^{kj} \omega_{\alpha, ij} (u)
{\pa^2 a \over \pa u^k \pa u^l} = \nn\\
&&
= \rho \, \eta^{kj} {\pa \omega_{\alpha, ij} \over \pa u^l}
 {\pa a \over \pa u^k} +
\rho \, \eta^{kj} \omega_{\alpha, ij} (u)
 \lambda \, \mu^{\gamma \beta}
\omega_{\gamma, kl} (u)  b_{\beta} (u) = \nn\\
&&
= \rho \, \eta^{kj} {\pa \omega_{\alpha, lj} \over \pa u^i}
 {\pa a \over \pa u^k} +
\rho \, \eta^{kj} \omega_{\alpha, lj} (u)
 \lambda \, \mu^{\gamma \beta}
\omega_{\gamma, ki} (u)  b_{\beta} (u),
\eea
whence we have
\be
{\pa \omega_{\alpha, ij} \over \pa u^l} =
{\pa \omega_{\alpha, lj} \over \pa u^i} \label{a2}
\ee
and
\be
\eta^{kj} \omega_{\alpha, ij} (u)
\omega_{\gamma, kl} (u) = \eta^{kj} \omega_{\alpha, lj} (u)
\omega_{\gamma, ki} (u). \label{b2}
\ee
It follows from (\ref{a1}) and (\ref{a2}) that
there locally exist some functions $\psi_{\alpha} (u),$
$1 \leq \alpha \leq L,$ such that
\be
\omega_{\alpha, ij} (u) = {\pa^2 \psi_{\alpha} \over \pa u^i \pa u^j}.
\ee
Then the relations (\ref{b1}) and (\ref{b2}) are
equivalent to the nonlinear system (\ref{g3}), (\ref{r3}) for the
functions $\psi_{\alpha} (u),$ $1 \leq \alpha \leq L.$

\bt  {\rm \cite{4m}}
The nonlinear system {\rm (\ref{g3}), (\ref{r3})} is
integrable by the inverse scattering method.
\et

\subsection{Integrable invariant description of flat torsionless
submanifolds in \\pseudo-Euclidean spaces}

In arbitrary local coordinates, we obtain the following
integrable description of all $N$-dimensional flat torsionless submanifolds
in $(N + L)$-dimensional pseudo-Euclidean spaces.

\bt {\rm \cite{3m}, \cite{4m}}
For each $N$-dimensional flat torsionless submanifold
in an $(N + L)$-dimensional
pseudo-Euclidean space with a flat first fundamental form $g_{ij} (u)$,
there locally
exist functions $\psi_{\alpha} (u),$ $1 \leq \alpha \leq L,$ such that
the second fundamental forms have the form
\be
(\omega_{\alpha})_{ij} (u) = \nabla_i \nabla_j \psi_{\alpha}, \label{f1}
\ee
where $\nabla_i$ is the
covariant differentiation defined by the Levi-Civita connection
generated by the metric $g_{ij} (u)$.
The class of $N$-dimensional flat torsionless submanifolds
in $(N + L)$-dimensional
pseudo-Euclidean spaces is described by the following
integrable system of nonlinear equations for
the functions $\psi_{\alpha} (u),$ $1 \leq \alpha \leq L${\rm :}
\be
 \sum_{n=1}^N  \nabla^n \nabla_i \psi_{\alpha}
\nabla_n \nabla_l \psi_{\beta}
=  \sum_{n=1}^N  \nabla^n \nabla_i \psi_{\beta}
\nabla_n \nabla_l \psi_{\alpha}, \label{non1}
\ee
\be
\sum_{\alpha = 1}^L \sum_{\beta = 1}^L \mu^{\alpha \beta}
\nabla_i \nabla_j \psi_{\alpha} \nabla_k \nabla_l \psi_{\beta} =
\sum_{\alpha = 1}^L \sum_{\beta = 1}^L \mu^{\alpha \beta}
\nabla_i \nabla_k
\psi_{\alpha} \nabla_j \nabla_l \psi_{\beta}, \label{non2}
\ee
where $\nabla_i$ is the
covariant differentiation defined by the Levi-Civita connection
generated by a flat metric
$g_{ij} (u)$, $\nabla^i = g^{is} (u) \nabla_s,$
$g^{is} (u) g_{sj} (u) = \delta^i_j.$
Moreover, in this case,
the systems
of hydrodynamic type
\be
u^i_{t_{\alpha}} = \left (\nabla^i \nabla_j
\psi_{\alpha}\right ) u^j_x,
\ \ \ \ 1 \leq \alpha \leq L, \label{h1}
\ee
are commuting integrable bi-Hamiltonian systems
of hydrodynamic type.

Any solution $\psi_{\alpha} (u),$ $1 \leq \alpha \leq L,$
of the integrable nonlinear system {\rm (\ref{non1}), (\ref{non2})}
determines a unique {\rm (}up to motions{\rm )}
$N$-dimensional flat torsionless
submanifold of
the corresponding $(N + L)$-dimensional pseudo-Euclidean space with the
first fundamental form $g_{ij} (u) d u^i d u^j$ and the
second fundamental forms {\rm (\ref{f1})}.
\et

\section{Reduction to the associativity equations of
two-dimensional \\topological quantum field theories
and potential flat \\torsionless submanifolds
in pseudo-Euclidean spaces} \label{section3}

\subsection{Special case of flat torsionless
submanifolds, when the Gauss and the Ricci
\\equations coincide, and the associativity equations}

We now also find some natural and very important
integrable reductions of the
nonlinear system (\ref{g3}), (\ref{r3}).
We show that the class of flat torsionless
submanifolds in pseudo-Euclidean spaces
is quite rich, and we describe a nontrivial and very
important family of submanifolds of this class.
This family is generated by the associativity
equations of two-dimensional topological quantum field theories
(the WDVV equations).
First of all, we note that although the Gauss equations {\rm (\ref{g3})}
and the Ricci equations {\rm (\ref{r3})} for
flat torsionless submanifolds in pseudo-Euclidean spaces
are essentially different, they are fantastically similar.
The case of a natural reduction under which
the Gauss equations {\rm (\ref{g3})} and the Ricci equations
{\rm (\ref{r3})} merely coincide is of particular interest.
Such a reduction readily leads to the associativity
equations of two-dimensional topological quantum field theories.

\bt  {\rm \cite{2m}--\cite{4m}}
If $L = N$, $\mu^{ij} = c \eta^{ij},$ $1 \leq i, j \leq N,$
$c$ is an arbitrary nonzero constant,
and $\psi_{\alpha} (u) = {\pa \Phi / \pa u^{\alpha}},$
$1 \leq \alpha \leq N,$ where $\Phi = \Phi (u^1, \ldots, u^N)$,
then the Gauss equations {\rm (\ref{g3})} coincide with the Ricci equations
{\rm (\ref{r3})}, and each of them coincides with the associativity
equations {\rm (\ref{ass1a})} of two-dimensional
topological quantum field theories
{\rm (}the WDVV equations{\rm )}
for the potential $\Phi (u)$.
\et

\subsection{Potential flat torsionless submanifolds
in pseudo-Euclidean spaces and \\the associativity equations}

\bd {\rm \cite{4m}
A flat torsionless $N$-dimensional submanifold in
a $2N$-dimensional pseudo-Euclidean space with a flat
first fundamental form $g_{ij} (u) d u^i d u^j$
is called {\it potential\/}
if there always locally exist a certain
function $\Phi (u)$ in a neighborhood on the submanifold such that
the second fundamental forms of this submanifold
locally in this neighborhood have the form
\be
(\omega_i)_{jk} (u) d u^j d u^k =
\left ( \nabla_i \nabla_j \nabla_k \Phi (u) \right ) d u^j d u^k,
\ \ \ \ 1 \leq i \leq N,
\ee
where $\nabla_i$ is the
covariant differentiation defined by the Levi-Civita connection
generated by the flat metric
$g_{ij} (u)$.
}
\ed

\bt  {\rm \cite{2m}--\cite{4m}}
The associativity
equations of two-dimensional topological quantum field theories
describe a special class of $N$-dimensional
flat submanifolds without torsion
in $2N$-dimensional
pseudo-Euclidean spaces, namely, exactly the class of {\rm potential flat
torsionless submanifolds}.
\et

According to the Bonnet theorem,
any solution $\Phi (u)$
of the associativity equations (\ref{ass1a})
(with an arbitrary fixed constant metric $\eta_{ij}$)
determines a unique (up to motions)
$N$-dimensional potential flat torsionless
submanifold of
the corresponding $2N$-dimensional pseudo-Euclidean space with the
first fundamental form $\eta_{ij} d u^i d u^j$ and the
second fundamental forms
$\omega_n (u) =
(\partial^3 \Phi / (\partial u^n \partial u^i \partial u^j))
d u^i d u^j$
 given by the third derivatives of the
potential $\Phi (u)$.
Here, we do not distinguish solutions of
the associativity equations (\ref{ass1a}) up to
terms quadratic in the coordinates $u$.

\bt {\rm \cite{3m}, \cite{4m}}
On each potential flat torsionless submanifold
in a pseudo-Euclidean space, there is a structure of a Frobenius
algebra given {\rm (}in flat coordinates{\rm )}
by the flat first fundamental form
$\eta_{ij}$ and the Weingarten operators $(A_s)^i_j (u) = -
\eta^{ik} (\omega_s)_{kj} (u)${\rm :}
\bea
&&
\langle e_i, e_j \rangle =
\eta_{ij}, \ \ \ \ e_i \circ e_j = c^k_{ij} (u) e_k,  \ \ \ \
e_i = {\pa \over \pa u^i}, \nn\\
&&
c^k_{ij} (u^1, \ldots, u^N) = - (A_i)^k_j (u) = \eta^{ks}
(\omega_i)_{sj} (u^1, \ldots, u^N).
\eea
In arbitrary local coordinates, this Frobenius structure has the form
\bea
&&
\langle e_i, e_j \rangle =
g_{ij}, \ \ \ \ e_i \circ e_j = c^k_{ij} (u) e_k,  \ \ \ \
e_i = {\pa \over \pa u^i}, \nn\\
&&
c^k_{ij} (u^1, \ldots, u^N) = - (A_i)^k_j (u) =
g^{ks} (u^1, \ldots, u^N)
(\omega_i)_{sj} (u^1, \ldots, u^N),
\eea
where $g^{ij} (u)$ is the contravariant metric inverse to the
first fundamental form $g_{ij} (u)$,
$g^{is} (u) g_{sj} (u) = \delta^i_j,$ and
$(\omega_k)_{ij} (u) d u^i d u^j,$ $1 \leq k \leq N,$ are
the second fundamental forms.
\et

\bt  {\rm \cite{3m}, \cite{4m}}
Each $N$-dimensional Frobenius manifold
can be locally represented as a potential flat torsionless
$N$-dimensional
submanifold in a $2N$-dimensional
pseudo-Euclidean space.
\et

\section{Realization of Frobenius manifolds as
submanifolds in \\pseudo-Euclidean
spaces}

It is important to note that we have at least two essentially
different possibilities for signature of the corresponding
ambient $2N$-dimensional pseudo-Euclidean space, namely,
we can always consider the ambient $2N$-dimensional
pseudo-Euclidean space of zero signature, and we can also consider
the ambient $2N$-dimensional
pseudo-Euclidean space whose signature
is equal to doubled
signature of the metric $\eta_{ij}$. Thus,
if the metric $\eta_{ij}$ of a Frobenius manifold has a nonzero
signature, then according to our construction
we have two essentially different
possibilities for realization
of the Frobenius manifold as a potential flat torsionless
submanifold.

\bt {\rm \cite{dga}}
For an arbitrary Frobenius manifold, which is locally
given by a solution
$\Phi (u^1, \ldots, u^N)$ of
the associativity equations {\rm (\ref{ass1a})}, the corresponding
potential flat torsionless submanifold in a $2N$-di\-mensional
pseudo-Euclidean space that realizes this Frobenius manifold
is given by the $2N$-component vector function $r (u^1, \ldots, u^N)$
satisfying the following compatible linear system of second-order
partial differential equations{\rm :}
\be
{\partial^2 r \over \partial u^i \partial u^j}
= c \eta^{kl}
{\partial^3 \Phi \over \partial u^i \partial u^j \partial u^k}
{\partial n \over \partial u^l},  \label{rf1}
\ee
\be
{\partial^2 n \over \partial u^i \partial u^j}
= - \eta^{kl}
{\partial^3 \Phi \over \partial u^i \partial u^j \partial u^k}
{\partial r \over \partial u^l},  \label{rf2}
\ee
where $n (u^1, \ldots, u^N)$ is a $2N$-component vector function,
$c$ is an arbitrary nonzero constant {\rm (}a deformation parameter
preserving the corresponding Frobenius structure{\rm )}. In particular,
two essentially different cases $c = 1$ and $c = - 1$ correspond
to ambient $2N$-dimensional
pseudo-Euclidean spaces of different signatures
{\rm (}if the metric $\eta_{ij}$ has a nonzero
signature{\rm )}. The consistency of the linear system
{\rm (\ref{rf1}), (\ref{rf2})} is equivalent to the associativity
equations {\rm (\ref{ass1a})}.
\et

\section{General nonlocal Hamiltonian operators
of hydrodynamic type}

Now we consider applications of our construction
to the theory of integrable systems, the theory of
nonlocal Hamiltonian operators of
hydrodynamic type, the theory of compatible Poisson structures
and the theory of bi-Hamiltonian integrable hierarchies of
hydrodynamic type. Recall that
the theory of nonlocal Hamiltonian operators of
hydrodynamic type was invented by the author and Ferapontov
in 1990--1991 \cite{2a}, \cite{1f}
in connection with vital necessities of
the Hamiltonian theory of systems of
hydrodynamic type proposed by Dubrovin and Novikov \cite{2dn}
and developed by Tsarev \cite{tsar}.
In this paper we give an integrable description
of all nonlocal Hamiltonian operators of hydrodynamic type
with flat metrics. This
nontrivial special class of Hamiltonian operators
is generated by flat torsionless submanifolds
in pseudo-Euclidean spaces and
closely connected with the associativity equations of
two-dimensional topological quantum field theories
and the theory of Frobenius manifolds.
The Hamiltonian operators of this class
are of special interest for many other reasons too.
In particular, any such
Hamiltonian operator always determines integrable
structural flows (some systems of hydrodynamic type),
always gives a nontrivial pencil of compatible
Hamiltonian operators and generates bi-Hamiltonian integrable
hierarchies of hydrodynamic type.
The affinors of any such
Hamiltonian operator generate some special integrals
in involution. The nonlinear systems describing
integrals in involution are of independent
great interest. The equations of associativity
of two-dimensional topological quantum field theories
(the WDVV equations) describe an important special
class of integrals in involution, a special class
of nonlocal Hamiltonian operators of hydrodynamic type
with flat metrics, a special class of
compatible local and nonlocal Poisson structures and important
special classes of bi-Hamiltonian integrable hierarchies of
systems of hydrodynamic type.
Moreover, we show that each flat torsionless submanifold
in a pseudo-Euclidean space (recall that this class of submanifolds
is described in our paper by an
integrable system \cite{4m}) gives a set of integrals in involution,
nontrivial local and nonlocal Hamiltonian operators of hydrodynamic type
with flat metrics, a pencil of compatible
Poisson structures and generates bi-Hamiltonian integrable hierarchies of
systems of hydrodynamic type.

Recall that general
nonlocal Hamiltonian operators of
hydrodynamic type, namely, Hamiltonian operators
of the form
\be
P^{ij} = g^{ij}(u(x)) {d \over dx} +
b^{ij}_k (u(x))\, u^k_x + \sum_{n =1}^L
\varepsilon^n (w_n)^i_k (u (x))
u^k_x \left ( {d \over dx} \right )^{-1} \circ
(w_n)^j_s (u (x)) u^s_x, \label{nonle}
\ee
where $\det (g^{ij} (u)) \neq 0,$
$\varepsilon^n = \pm 1,$ $1 \leq n \leq L,$ $u^1, \ldots, u^N$
are local coordinates, $u = (u^1, \ldots, u^N),$
$u^i (x),$ $1 \leq i \leq N,$ are functions (fields)
of one independent variable $x$, and the coefficients
$g^{ij} (u),$ $b^{ij}_k (u),$ $(w_n)^i_j (u),$
$1 \leq i, j, k \leq N,$ $1 \leq n \leq L,$ are smooth
functions of local coordinates,
were studied by Ferapontov in \cite{1f}
(see also \cite{2dn}, \cite{2a}).

Hamiltonian operators of the general form
(\ref{nonle}) (local and nonlocal) play a key role
in the Hamiltonian theory of systems of hydrodynamic type
\cite{2dn}--\cite{1f}.
Recall that an operator $M^{ij}$ is said to be
{\it Hamiltonian\/} if the operator defines a Poisson bracket
\be
\{ I, J \} = \int {\delta I \over \delta u^i (x)} M^{ij}
{\delta J \over \delta u^j (x)} dx \label{brac}
\ee
on arbitrary functionals $I$ and $J$ on the space of the fields
$u^i (x)$, i.e., the bracket
(\ref{brac}) is skew-symmetric and
satisfies the Jacobi identity.

It was proved in \cite{1f} that
the operator (\ref{nonle}) is Hamiltonian if and only if
$g^{ij} (u)$ is a symmetric (pseudo-Riemannian)
contravariant metric and the following relations are satisfied for
the coefficients of the operator:

1) $b^{ij}_k (u) = - g^{is} (u) \Gamma ^j_{sk} (u),$
where
$\Gamma^j_{sk} (u)$ is the Levi-Civita connection
generated by the contravariant metric
$g^{ij} (u)$,

2) $g^{ik} (u) (w_n)^j_k (u) =
g^{jk} (u) (w_n)^i_k (u),$

3) $\nabla_k (w_n)^i_j (u) =
\nabla_j (w_n)^i_k (u),$ where $\nabla_k$
is the covariant differentiation
generated by the Levi-Civita connection
$\Gamma^j_{sk} (u)$ of the metric $g^{ij} (u),$

4) $R^{ij}_{kl} (u) =  \sum_{n =1}^L
\varepsilon^n
\left ( (w_n)^i_l (u) (w_n)^j_k (u)
- (w_n)^j_l (u) (w_n)^i_k (u) \right ),$ where
$$R^{ij}_{kl} (u)= g^{is} (u) R^j_{skl} (u)$$
is the Riemannian curvature tensor of the metric
$g^{ij} (u),$

5) $[w_n (u), w_m (u)] = 0$, i.e.,
the family $(w_n)^i_j (u),$ $1 \leq n \leq L,$ of (1, 1)-tensors
({\it affinors\/})
is commutative.

Each Hamiltonian operator of the form
(\ref{nonle}) exactly corresponds to an $N$-dimensional
submanifold with flat normal bundle embedded in
a pseudo-Euclidean space $E^{N + L}$.
Here, the covariant metric $g_{ij} (u)$ (for which
$g_{is} (u) g^{sj} (u) = \delta^j_i$) is the first
fundamental form, and the affinors
$w_n (u),$ $1 \leq n \leq L,$ are the
Weingarten operators of this embedded submanifold
($g_{is} (u) (w_n)^s_j (u)$ are the corresponding
second fundamental forms).
Correspondingly, the relations 2)--4)
are the Gauss--Codazzi
equations for an $N$-dimensional submanifold
with zero torsion embedded in a pseudo-Euclidean
space $E^{N + L}$ \cite{1f}. The relations 5)
are equivalent to the Ricci equations for
this embedded submanifold.

Having in mind further applications
to arbitrary Frobenius manifolds, we prefer
to consider general nonlocal Hamiltonian operators
of hydrodynamic type in the form
\be
P^{ij} = g^{ij}(u(x)) {d \over dx} +
b^{ij}_k (u(x))\, u^k_x + \sum_{m =1}^L \sum_{n =1}^L
\mu^{mn} (w_m)^i_k (u (x))
u^k_x \left ( {d \over dx} \right )^{-1} \circ
(w_n)^j_s (u (x)) u^s_x, \label{nonlm}
\ee
where $\det (g^{ij} (u)) \neq 0,$
$\mu^{mn}$ is an arbitrary nondegenerate symmetric
constant matrix.
Each operator of the form (\ref{nonlm}) can be
reduced
to the form (\ref{nonle}) (and conversely,
each operator of the form (\ref{nonlm}) can be
obtained from some operator of the form (\ref{nonle}))
by a linear transformation
$w_n (u) = c^l_n \widetilde w_l (u)$
in the
vector space of affinors $w_n (u),$ $1 \leq n \leq L$;
here $c^l_n$ is
an arbitrary nondegenerate constant matrix.
Among all the conditions 1)--5) for the Hamiltonian property
of the operator
(\ref{nonle}), these transformations change only the condition
4) for the Riemannian curvature tensor of the metric.
The condition 4) for the operator
(\ref{nonlm}) takes the form
$$R^{ij}_{kl} (u) =  \sum_{m =1}^L \sum_{n =1}^L
\mu^{mn}
\left ( (w_m)^i_l (u) (w_n)^j_k (u)
- (w_m)^j_l (u) (w_n)^i_k (u) \right ),$$
and all the other conditions 1)--3) and 5) for the Hamiltonian property
remain unchanged.

Consider all the relations for the coefficients
of the nonlocal Hamiltonian operator
(\ref{nonlm}) in a form convenient for further use.

\bl {\rm \cite{3m}}
The operator {\rm (\ref{nonlm})}, where $\det (g^{ij} (u)) \neq 0,$
is Hamiltonian if and only
if its coefficients satisfy the relations
\be
g^{ij} = g^{ji}, \label{01}
\ee
\be
{\pa g^{ij} \over \pa u^k} = b^{ij}_k + b^{ji}_k, \label{02}
\ee
\be
g^{is} b^{jk}_s = g^{js} b^{ik}_s, \label{03}
\ee
\be
g^{is} (w_n)^j_s = g^{js} (w_n)^i_s, \label{04}
\ee
\be
(w_n)^i_s (w_m)^s_j =
(w_m)^i_s (w_n)^s_j, \label{05}
\ee
\be
g^{is} g^{jr} {\pa (w_n)^k_r \over \pa u^s}
- g^{jr} b^{ik}_s (w_n)^s_r  =
g^{js} g^{ir} {\pa (w_n)^k_r \over \pa u^s}
- g^{ir} b^{jk}_s (w_n)^s_r, \label{06}
\ee
\be
g^{is} \left ( {\pa b^{jk}_s \over \pa u^r}
- {\pa b^{jk}_r \over \pa u^s} \right )
+ b^{ij}_s b^{sk}_r - b^{ik}_s b^{sj}_r =
\sum_{m = 1}^L \sum_{n =1}^L \mu^{mn} g^{is}
\left ( (w_m)^j_r (w_n)^k_s -
 (w_m)^j_s (w_n)^k_r \right ). \label{07}
\ee
\el

\section{Nonlocal Hamiltonian operators of hydrodynamic type
with flat \\metrics and special pencils
of Hamiltonian operators} \label{pencil}

Let us consider the important special case
of the nonlocal Hamiltonian operators of the form
(\ref{nonlm}) when the metric $g^{ij} (u)$ is flat.
Recall that each flat metric uniquely determines
a local Hamiltonian operator of hydrodynamic type
(i.e., a Hamiltonian operator of the form (\ref{nonlm})
with zero affinors) known as
a Dubrovin--Novikov Hamiltonian operator \cite{2dn}.
We prove that for each flat metric
there also exist a remarkable class of nonlocal
Hamiltonian operators of hydrodynamic type
with this flat metric and nontrivial affinors;
moreover, these Hamiltonian operators
have important applications in the theory of
Frobenius manifolds and integrable hierarchies.
First of all, note the following important
property of nonlocal Hamiltonian operators
of hydrodynamic type with flat metrics.
Recall that two Hamiltonian operators
are said to be {\it compatible\/} if any linear combination
of these Hamiltonian operators is also a Hamiltonian operator \cite{5},
i.e., they form a {\it pencil of Hamiltonian operators\/}
and, correspondingly, they form a {\it pencil of
Poisson brackets}.

\bl {\rm \cite{3m}} \label{fl}
The metric $g^{ij} (u)$ of a Hamiltonian operator
of the form {\rm (\ref{nonlm})} is flat if and only if
this operator defines the pencil
\bea
P^{ij}_{\lambda_1, \lambda_2} &=&
\lambda_1 \left ( g^{ij}(u(x)) {d \over dx} +
b^{ij}_k (u(x))\, u^k_x \right ) + \nn\\
&+&
\lambda_2 \sum_{m =1}^L \sum_{n =1}^L
\mu^{mn} (w_m)^i_k (u (x))
u^k_x \left ( {d \over dx} \right )^{-1} \circ
(w_n)^j_s (u (x)) u^s_x, \label{nonlma}
\eea
of compatible Hamiltonian operators,
where $\lambda_1$ and $\lambda_2$ are arbitrary constants.
\el

Indeed, if the operator (\ref{nonlm})
is Hamiltonian, then its coefficients satisfy the relations
(\ref{01})--(\ref{07}). It is obvious that in this case
the relations (\ref{01})--(\ref{06}) for the operator
(\ref{nonlma}) are always satisfied for any constants
$\lambda_1$ and $\lambda_2$, and the relation (\ref{07})
is satisfied for any constants $\lambda_1$ and $\lambda_2$
if and only if the left- and right-hand sides of this relation
are zero identically.

It follows from the relations (\ref{01})--(\ref{03})
for the Hamiltonian operator
(\ref{nonlm}) that the Riemannian curvature tensor
of the metric $g^{ij} (u)$ has the form
\be
R^{ijk}_r (u) = g^{is} (u) R^{jk}_{sr} (u) =
g^{is} (u) \left ( {\pa b^{jk}_s \over \pa u^r}
- {\pa b^{jk}_r \over \pa u^s} \right )
+ b^{ij}_s (u) b^{sk}_r (u) - b^{ik}_s (u) b^{sj}_r (u). \label{kr}
\ee
Consequently, if the metric $g^{ij} (u)$
of a Hamiltonian operator of the form (\ref{nonlm}) is flat,
i.e., $R^{ijk}_r (u) = 0,$ then the relation (\ref{07})
becomes
$$\sum_{m = 1}^L \sum_{n =1}^L \mu^{mn} g^{is}
\left ( (w_m)^j_r (u) (w_n)^k_s (u) -
 (w_m)^j_s (u) (w_n)^k_r (u) \right ) = 0.$$
Thus the metric $g^{ij} (u)$ of a Hamiltonian operator
of the form (\ref{nonlm}) is flat if and only if
the left- and right-hand sides of the relation
(\ref{07}) for the Hamiltonian operator
(\ref{nonlm}) are zero identically.
In this case, the left-
and right-hand sides of the relation
(\ref{07}) for the operator
(\ref{nonlma}) are also zero identically
for any constants $\lambda_1$ and $\lambda_2$,
i.e., we obtain a pencil of compatible
Hamiltonian operators (\ref{nonlma}).
Note also that for the pencil of Hamiltonian
operators $P^{ij}_{\lambda_1, \lambda_2}$ given by
the formula (\ref{nonlma})
it readily follows from the Dubrovin--Novikov
theorem \cite{2dn} applied to the local operator
$P^{ij}_{1, 0}$ that the metric
$g^{ij} (u)$ is flat. Lemma \ref{fl} is proved.

Thus if the metric
$g^{ij} (u)$ of a Hamiltonian operator
of the form (\ref{nonlm}) is flat, then the operator
\be
P^{ij}_{0, 1} =
 \sum_{m =1}^L \sum_{n =1}^L
\mu^{mn} (w_m)^i_k (u (x))
u^k_x \left ( {d \over dx} \right )^{-1} \circ
(w_n)^j_s (u (x)) u^s_x \label{nonlmb}
\ee
is also a Hamiltonian operator obtained by the degeneration
as $\lambda_1 \rightarrow 0$. Moreover,
in this case, this Hamiltonian operator
is always compatible with the local Hamiltonian
operator of hydrodynamic type
(the Dubrovin--Novikov operator)
\be
P^{ij}_{1, 0} =
g^{ij}(u(x)) {d \over dx} +
b^{ij}_k (u(x))\, u^k_x.  \label{nonlmc}
\ee

The compatible Hamiltonian operators
(\ref{nonlmb}) and
(\ref{nonlmc}) always generate the corresponding
integrable bi-Hamiltonian hierarchies. We construct
these integrable hierarchies further in Section \ref{hier5}.

\section{Integrability of structural flows} \label{potoki}

We recall that systems of hydrodynamic type
\be
u^i_{t_n} = (w_n)^i_j (u) u^j_x, \ \ \ \ 1 \leq n \leq L, \label{strhy}
\ee
are called {\it structural flows} of the nonlocal
Hamiltonian operator of hydrodynamic type (\ref{nonlm})
(see \cite{1f}, \cite{4}).

\bl  {\rm \cite{3m}}
All the structural flows {\rm (\ref{strhy})} of
any nonlocal Hamiltonian operator
of hydrodynamic type with flat metric
are commuting integrable bi-Hamiltonian
systems of hydrodynamic type.
\el

Maltsev and Novikov proved in \cite{4}
(see also \cite{1f}) that the structural flows of any
nonlocal Hamiltonian operator of
hydrodynamic type (\ref{nonlm}) are always Hamiltonian
with respect to this Hamiltonian operator.
Let us consider an arbitrary nonlocal
Hamiltonian operator of hydrodynamic type
(\ref{nonlm}) with a flat metric $g^{ij} (u)$
and
the pencil of compatible Hamiltonian operators
(\ref{nonlma}) corresponding to this Hamiltonian
operator. The corresponding structural flows
are necessarily Hamiltonian with respect to each of the
operators in the Hamiltonian pencil
(\ref{nonlma}) and, consequently,
they are integrable bi-Hamiltonian systems.

\section{Integrable description of nonlocal Hamiltonian
operators of \\hydrodynamic type with
flat metrics} \label{nonop}

Let us describe all the nonlocal
Hamiltonian operators of hydrodynamic type
with flat metrics.
The form of the Hamiltonian operator (\ref{nonlm})
is invariant with respect to local changes of coordinates,
and also all the coefficients of the operator are
transformed as the corresponding differential-geometric objects.
Since the metric is flat, there exist local coordinates
in which the metric is reduced to a constant matrix
$\eta^{ij},$ $\eta^{ij}= {\rm const},$ $\det (\eta^{ij}) \neq 0,$
$\eta^{ij} = \eta^{ji}$.
In these local coordinates, all the coefficients of
the Levi-Civita connection are zero, and
the Hamiltonian operator has the form
\be
\widetilde P^{ij} = \eta^{ij} {d \over dx} +
\sum_{m =1}^L \sum_{n =1}^L
\mu^{mn} (\widetilde w_m)^i_k (u (x))
u^k_x \left ( {d \over dx} \right )^{-1} \circ
(\widetilde w_n)^j_s (u (x)) u^s_x. \label{nonl3}
\ee
Description of nonlocal Hamiltonian operators of hydrodynamic type with
flat metrics coincides with description of flat torsionless
submanifolds in pseudo-Euclidean spaces.

\bt  {\rm \cite{2m}, \cite{3m}}
The operator {\rm (\ref{nonl3})}, where $\eta^{ij}$ and $\mu^{mn}$
are arbitrary nondegenerate symmetric constant matrices,
is Hamiltonian if and only if
there exist functions $\psi_n (u),$ $1 \leq n \leq L,$
such that
\be
(\widetilde w_n)^i_j (u) = \eta^{is}
{\pa^2 \psi_n \over \pa u^s \pa u^j}, \label{n}
\ee
and the integrable system {\rm (\ref{g3}), (\ref{r3})}
of nonlinear equations describing all flat torsionless
submanifolds in pseudo-Euclidean spaces
is satisfied.
\et

The relations (\ref{01})--(\ref{03}) for any operator
of the form (\ref{nonl3}) are automatically fulfilled,
and the relation (\ref{06}) for any operator of the form
(\ref{nonl3}) has the form
\be
{\pa (\widetilde w_n)^k_r \over \pa u^s}
 =
{\pa (\widetilde w_n)^k_s \over \pa u^r}, \label{06a}
\ee
and, consequently, there locally exist functions
$\varphi_n^i (u),$ $1 \leq i \leq N,$ $1 \leq n \leq L,$
such that
\be
(\widetilde w_n)^i_j (u) = {\pa \varphi_n^i \over \pa u^j}.
\ee
Then relation (\ref{04}) becomes
\be
\eta^{is} {\pa \varphi_n^j \over \pa u^s} =
\eta^{js} {\pa \varphi_n^i \over \pa u^s}
\ee
or, equivalently,
\be
{\pa \left ( \eta_{is} \varphi_n^s \right ) \over \pa u^j} =
{\pa \left ( \eta_{js} \varphi_n^s \right ) \over \pa u^i}, \label{04a}
\ee
where the matrix $\eta_{ij}$ is inverse to the matrix
$\eta^{ij}$,
$\eta_{is} \eta^{sj} = \delta^j_i.$
It follows from the relation (\ref{04a}) that there locally
exist functions $\psi_n (u),$ $1 \leq n \leq L,$ such that
\be
\eta_{is} \varphi_n^s = {\pa \psi_n \over \pa u^i}.
\ee
Thus
\be
\varphi_n^i = \eta^{is} {\pa \psi_n \over \pa u^s},
\ \ \ (\widetilde w_n)^i_j (u) = \eta^{is}
{\pa^2 \psi_n \over \pa u^s \pa u^j}.
\ee
In this case, the relations (\ref{05}) and (\ref{07})
become  (\ref{r3}) and (\ref{g3}) respectively.

The nonlinear equations (\ref{g3}) and (\ref{r3})
describing all nonlocal Hamiltonian operators of
hydrodynamic type with flat metrics are exactly equivalent
to the conditions that a flat $N$-dimensional submanifold
with flat normal bundle, with the first fundamental form
$\eta_{ij} d u^i d u^j$ and the second fundamental
forms $\omega_n (u)$ given by Hessians of $L$ functions
$\psi_n (u),$ $1 \leq n \leq L$,
$$
\omega_n (u) = {\pa^2 \psi_n \over \pa u^i \pa u^j} d u^i d u^j,
$$
is embedded in an $(N + L)$-dimensional pseudo-Euclidean space.

\section{Integrable description of a special class of
pencils of \\Hamiltonian
operators} \label{intpen}

Now we give an integrable description of a special class of
pencils of Hamiltonian operators and a special class of
integrable bi-Hamiltonian hierarchies of hydrodynamic type.

\bt  {\rm \cite{3m}}
If functions $\psi_n (u),$ $1 \leq n \leq L,$
are a solution of the integrable nonlinear system
{\rm (\ref{g3}), (\ref{r3})}, then the systems
of hydrodynamic type {\rm (}the structural flows of the
corresponding nonlocal Hamiltonian operator of hydrodynamic type with
flat metric{\rm )}
\be
u^i_{t_n} = \eta^{is}
{\pa^2 \psi_n \over \pa u^s \pa u^j} u^j_x,
\ \ \ \ 1 \leq n \leq L, \label{str2hy}
\ee
are commuting integrable bi-Hamiltonian systems
of hydrodynamic type. Moreover, in this case
the nonlocal operator
\be
M^{ij}_1 =
 \sum_{m =1}^L \sum_{n =1}^L
\mu^{mn} \eta^{ip} \eta^{jr} {\pa^2 \psi_m \over \pa u^p \pa u^k}
u^k_x \left ( {d \over dx} \right )^{-1} \circ
{\pa^2 \psi_n \over \pa u^r \pa u^s} u^s_x \label{nonlmd}
\ee
is also a Hamiltonian operator,
and this nonlocal Hamiltonian operator is compatible
with the constant Hamiltonian operator
\be
M^{ij}_2 =
\eta^{ij} {d \over dx}.  \label{nonlme}
\ee
\et

In arbitrary local coordinates, we obtain the following
integrable description of all nonlocal Hamiltonian
operators of hydrodynamic type with flat metrics and
the corresponding pencils of Hamiltonian operators.

\bt {\rm \cite{3m}}
The operator {\rm (\ref{nonlm})} with a flat metric $g^{ij} (u)$
is Hamiltonian if and only if
$b^{ij}_k (u) = - g^{is} (u) \Gamma^j_{sk} (u),$
where
$\Gamma^j_{sk} (u)$ is the flat connection
generated by the flat metric $g^{ij} (u)$, and there locally
exist functions $\psi_n (u),$ $1 \leq n \leq L,$ such that
\be
(w_n)^i_j (u) = \nabla^i \nabla_j \psi_n, \label{n2}
\ee
and the integrable system {\rm (\ref{non1}), (\ref{non2})} of
nonlinear equations describing all flat torsionless
submanifolds in pseudo-Euclidean spaces
is satisfied.
In particular, in this case the operator
\bea
M^{ij}_{\lambda_1, \lambda_2} &=&
\lambda_1 \left ( g^{ij}(u(x)) {d \over dx} -
g^{is} (u (x)) \Gamma^j_{sk} (u(x))\, u^k_x \right ) + \nn\\
&+&
\lambda_2
 \sum_{m =1}^L \sum_{n =1}^L
\mu^{mn} \nabla^i \nabla_k \psi_m
u^k_x \left ( {d \over dx} \right )^{-1} \circ
\nabla^j \nabla_s \psi_n u^s_x \label{nonlmdn}
\eea
is a Hamiltonian operator for any constants
$\lambda_1$ and $\lambda_2$, and the systems
of hydrodynamic type
\be
u^i_{t_n} = \nabla^i \nabla_j
\psi_n u^j_x,
\ \ \ \ 1 \leq n \leq L, \label{str2n}
\ee
are always commuting integrable bi-Hamiltonian systems
of hydrodynamic type.
\et

Hence, each flat torsionless submanifold in a pseudo-Euclidean space
generates a nonlocal Hamiltonian operator of hydrodynamic type with
flat metric, gives a special class of
pencils of Hamiltonian operators and special
integrable bi-Hamiltonian hierarchies of hydrodynamic type.
Now we construct an infinite
integrable bi-Hamiltonian hierarchy of hydrodynamic type
generated by an arbitrary
flat torsionless submanifold in a pseudo-Euclidean space.

\section{Integrable hierarchy generated by an arbitrary
flat torsionless \\submanifold in pseudo-Euclidean space} \label{hier5}

Consider the recursion operator
\be
R^i_j = \left ( M_1 (M_2)^{- 1} \right )^i_j =
\sum_{m =1}^L \sum_{n =1}^L
\mu^{mn} \eta^{ip} {\pa^2 \psi_m \over \pa u^p \pa u^k}
u^k_x \left ( {d \over dx} \right )^{-1} \circ
{\pa^2 \psi_n \over \pa u^j \pa u^s} u^s_x
\left ( {d \over dx} \right )^{-1} \label{rek}
\ee
corresponding to the compatible Hamiltonian operators
(\ref{nonlmd}) and (\ref{nonlme}).
Let us apply this recursion operator
(\ref{rek}) to the system of translations
with respect to $x$,
\be
u^i_t = u^i_x.
\ee
Each system in the hierarchy
\be
u^i_{t_s} = (R^s)^i_j u^j_x, \ \ \ \ s \in {\mathbb{Z}},  \label{ierar}
\ee
is a multi-Hamiltonian integrable system of
hydrodynamic type. In particular,
each system of the form
\be
u^i_{t_1} = R^i_j u^j_x,
\ee
i.e., the system
\be
u^i_{t_1} =
\sum_{m =1}^L \sum_{n =1}^L
\mu^{mn} \eta^{ip} {\pa^2 \psi_m \over \pa u^p \pa u^k}
u^k_x \left ( {d \over dx} \right )^{-1} \circ
{\pa^2 \psi_n \over \pa u^j \pa u^s} u^j u^s_x, \label{s1}
\ee
is integrable.

Since
\be
{\pa \over \pa u^r} \left ( {\pa^2 \psi_n \over
\pa u^j \pa u^s} u^j \right ) =
{\pa^3 \psi_n \over \pa u^j \pa u^s \pa u^r} u^j
+ {\pa^2 \psi_n \over \pa u^r \pa u^s} =
{\pa \over \pa u^s} \left ( {\pa^2 \psi_n \over
\pa u^j \pa u^r} u^j \right ),
\ee
there locally exist functions $F_n (u),$ $1 \leq n \leq L,$
such that
\be
{\pa^2 \psi_n \over
\pa u^j \pa u^s} u^j = {\pa F_n \over \pa u^s},
\ \ \ \ F_n = {\pa \psi_n \over \pa u^j} u^j - \psi_n.
\ee
Thus the system of hydrodynamic type
(\ref{s1}) has the local form
\be
u^i_{t_1} =
\sum_{m =1}^L \sum_{n =1}^L
\mu^{mn} \eta^{ip} F_n (u)
{\pa^2 \psi_m \over \pa u^p \pa u^k}
u^k_x. \label{s2}
\ee
This system of hydrodynamic type is
bi-Hamiltonian with respect to the compatible
Hamiltonian operators
(\ref{nonlmd}) and (\ref{nonlme}):
\be
u^i_{t_1} =
\sum_{m =1}^L \sum_{n =1}^L
\mu^{mn} \eta^{ip} \eta^{jr} {\pa^2 \psi_m \over \pa u^p \pa u^k}
u^k_x \left ( {d \over dx} \right )^{-1} \left (
{\pa^2 \psi_n \over \pa u^r \pa u^s} u^s_x {\delta H_1
\over \delta u^j (x)} \right ), \label{s3}
\ee
\be
H_1 = \int h_1 (u(x)) dx, \ \ \ \ h_1 (u(x)) =
{1 \over 2} \eta_{ij} u^i (x) u^j (x),
\ee
\be
u^i_{t_1} = \eta^{ij}  {d \over dx}  {\delta H_2
\over \delta u^j (x)}, \ \ \ \ H_2 = \int h_2 (u(x)) dx, \label{s4}
\ee
since in our case there always locally exists a function
$h_2 (u)$ such that
\be
\sum_{m =1}^L \sum_{n =1}^L
\mu^{mn} {\pa^2 \psi_m \over \pa u^j \pa u^k} F_n (u) =
{\pa^2 h_2 \over \pa u^j \pa u^k}.
\ee

Indeed, we have
\bea
&&
{\pa \over \pa u^i} \left ( \sum_{m =1}^L \sum_{n =1}^L
\mu^{mn} {\pa^2 \psi_m \over \pa u^j \pa u^k} F_n (u) \right ) =
\sum_{m =1}^L \sum_{n =1}^L
\mu^{mn} {\pa^3 \psi_m \over \pa u^i \pa u^j \pa u^k} F_n (u) + \nn\\
&&
+ \sum_{m =1}^L \sum_{n =1}^L
\mu^{mn} {\pa^2 \psi_m \over \pa u^j \pa u^k} {\pa F_n  \over
\pa u^i} = \sum_{m =1}^L \sum_{n =1}^L
\mu^{mn} {\pa^3 \psi_m \over \pa u^i \pa u^j \pa u^k} F_n (u)
+ \nn\\
&&
+ \sum_{m =1}^L \sum_{n =1}^L
\mu^{mn} {\pa^2 \psi_m \over \pa u^j \pa u^k} {\pa^2 \psi_n  \over
\pa u^i \pa u^s} u^s = \sum_{m =1}^L \sum_{n =1}^L
\mu^{mn} {\pa^3 \psi_m \over \pa u^i \pa u^j \pa u^k} F_n (u)
+ \nn\\
&&
+ \sum_{m =1}^L \sum_{n =1}^L
\mu^{mn} {\pa^2 \psi_m \over \pa u^j \pa u^i} {\pa^2 \psi_n  \over
\pa u^k \pa u^s} u^s, \label{k1}
\eea
where we have used the relation  (\ref{g3}).
Consequently, by virtue of symmetry
with respect to the indices $i$ and $j$, we obtain
\be
{\pa \over \pa u^i} \left ( \sum_{m =1}^L \sum_{n =1}^L
\mu^{mn} {\pa^2 \psi_m \over \pa u^j \pa u^k} F_n (u) \right ) =
{\pa \over \pa u^j} \left ( \sum_{m =1}^L \sum_{n =1}^L
\mu^{mn} {\pa^2 \psi_m \over \pa u^i \pa u^k} F_n (u) \right ),
\ee
i.e., there locally exist functions $a_k (u),$ $1 \leq k \leq N,$
such that
\be
 \sum_{m =1}^L \sum_{n =1}^L
\mu^{mn} {\pa^2 \psi_m \over \pa u^j \pa u^k} F_n (u) =
{\pa a_k \over \pa u^j}.
\ee
By virtue of symmetry with respect to
the indices $j$ and $k$, we obtain
\be
{\pa a_k \over \pa u^j} = {\pa a_j \over \pa u^k},
\ee
i.e., there locally exists a function
$h_2 (u)$ such that
\be
a_k (u) = {\pa h_2 \over \pa u^k}.
\ee
Thus
\be
 \sum_{m =1}^L \sum_{n =1}^L
\mu^{mn} {\pa^2 \psi_m \over \pa u^j \pa u^k} F_n (u) =
{\pa a_k \over \pa u^j} = {\pa^2 h_2 \over \pa u^j \pa u^k}. \label{k11}
\ee

Consider the next equation in the integrable hierarchy
(\ref{ierar}):
\bea
&&
u^i_{t_2} = (R^2)^i_j u^j_x = \nn\\
&&
= \sum_{m =1}^L \sum_{n =1}^L
\mu^{mn} \eta^{ip} {\pa^2 \psi_m \over \pa u^p \pa u^k}
u^k_x \left ( {d \over dx} \right )^{-1} \circ
{\pa^2 \psi_n \over \pa u^j \pa u^s} u^s_x
\left ( {d \over dx} \right )^{-1} \circ
\eta^{jr}  {d \over dx}  {\delta H_2
\over \delta u^r (x)} = \nn\\
&&
= \sum_{m =1}^L \sum_{n =1}^L
\mu^{mn} \eta^{ip} {\pa^2 \psi_m \over \pa u^p \pa u^k}
u^k_x \left ( {d \over dx} \right )^{-1} \circ
{\pa^2 \psi_n \over \pa u^j \pa u^s} u^s_x
\eta^{jr}  {\pa h_2
\over \pa u^r}.
 \label{ierar2}
\eea

Let us prove that in our case there always locally exist
functions $G_n (u),$ $1 \leq n \leq L,$ such that
\be
{\pa^2 \psi_n \over
\pa u^j \pa u^s} \eta^{jr}  {\pa h_2
\over \pa u^r} = {\pa G_n \over \pa u^s}. \label{gn}
\ee

Indeed, we have
\bea
&&
{\pa \over \pa u^p} \left ( {\pa^2 \psi_n \over
\pa u^j \pa u^s} \eta^{jr}  {\pa h_2
\over \pa u^r} \right ) = {\pa^3 \psi_n \over
\pa u^j \pa u^s \pa u^p} \eta^{jr}  {\pa h_2
\over \pa u^r}
+ {\pa^2 \psi_n \over
\pa u^j \pa u^s} \eta^{jr}  {\pa^2 h_2
\over \pa u^r \pa u^p} = \nn\\
&&
= {\pa^3 \psi_n \over
\pa u^j \pa u^s \pa u^p} \eta^{jr}  {\pa h_2
\over \pa u^r}
+ {\pa^2 \psi_n \over
\pa u^j \pa u^s} \eta^{jr} \left (\sum_{k =1}^L \sum_{l =1}^L
\mu^{kl} {\pa^2 \psi_k \over \pa u^r \pa u^p} F_l (u) \right )
= \nn\\
&&
= {\pa^3 \psi_n \over
\pa u^j \pa u^s \pa u^p} \eta^{jr}  {\pa h_2
\over \pa u^r}
+ \sum_{k =1}^L \sum_{l =1}^L
\mu^{kl} \eta^{jr} {\pa^2 \psi_n \over
\pa u^s \pa u^j} {\pa^2 \psi_k \over \pa u^r \pa u^p} F_l (u)
= \nn\\
&&
= {\pa^3 \psi_n \over
\pa u^j \pa u^s \pa u^p} \eta^{jr}  {\pa h_2
\over \pa u^r}
+ \sum_{k =1}^L \sum_{l =1}^L
\mu^{kl} \eta^{jr} {\pa^2 \psi_k \over
\pa u^s \pa u^j} {\pa^2 \psi_n \over \pa u^r \pa u^p} F_l (u)
= \nn\\
&&
= {\pa^3 \psi_n \over
\pa u^j \pa u^s \pa u^p} \eta^{jr}  {\pa h_2
\over \pa u^r}
+ \sum_{k =1}^L \sum_{l =1}^L
\mu^{kl} \eta^{jr} {\pa^2 \psi_k \over
\pa u^s \pa u^r} {\pa^2 \psi_n \over \pa u^j \pa u^p} F_l (u), \label{k2}
\eea
where we have used the relation (\ref{r3}) and
the symmetry of the matrix
$\eta^{jr}$. Thus we have proved that
the expression under consideration is
symmetric with respect to the indices
$p$ and $s$, i.e.,
\be
{\pa \over \pa u^p} \left ( {\pa^2 \psi_n \over
\pa u^j \pa u^s} \eta^{jr}  {\pa h_2
\over \pa u^r} \right ) =
{\pa \over \pa u^s} \left ( {\pa^2 \psi_n \over
\pa u^j \pa u^p} \eta^{jr}  {\pa h_2
\over \pa u^r} \right ). \label{k22}
\ee
Consequently, there locally exist functions
$G_n (u),$ $1 \leq n \leq L,$ such that
the relation
(\ref{gn}) is satisfied, and therefore, we have proved
that the second flow in the integrable hierarchy (\ref{ierar})
has the form of a local system of hydrodynamic type
\be
u^i_{t_2} = \sum_{m =1}^L \sum_{n =1}^L
\mu^{mn} \eta^{ip} G_n (u) {\pa^2 \psi_m \over \pa u^p \pa u^k}
u^k_x.
 \label{ierar2l}
\ee

Repeating the preceding argument word for word, we prove
by induction that if the functions
$\psi_n (u),$ $1 \leq n \leq L,$ are a solution of
the system of equations
 (\ref{g3}), (\ref{r3}), then for each $s \geq 1$
and for the corresponding function $h_s (u(x))$
(starting from the function
$h_1 (u(x)) = {1 \over 2} \eta_{ij} u^i (x) u^j (x)$)
there always locally exist functions
$F_n^{(s)} (u),$ $1 \leq n \leq L,$ such that
\be
{\pa^2 \psi_n \over
\pa u^j \pa u^p} \eta^{jr}  {\pa h_s
\over \pa u^r} = {\pa F_n^{(s)} \over \pa u^p} \label{fsn}
\ee
and there always locally exists a function $h_{s+1} (u(x))$
such that
\be
\sum_{m =1}^L \sum_{n =1}^L
\mu^{mn} {\pa^2 \psi_m \over \pa u^j \pa u^k} F_n^{(s)} (u) =
{\pa^2 h_{s+1} \over \pa u^j \pa u^k}.   \label{hsn}
\ee
Above we have already proved that this statement
is true for $s = 1$
(in this case, in particular, $F_n^{(1)} = F_n,$ $F_n^{(2)} = G_n$).
It can be proved in just the same way that
if this statement is true for $s = K \geq 1$,
then it is true also for $s = K + 1$
(see (\ref{k1})--(\ref{k11}) and (\ref{k2}), (\ref{k22})).
Thus we have proved that for each $s \geq 1$
the corresponding flow of the integrable hierarchy (\ref{ierar})
has the form of a local system of hydrodynamic type
\be
u^i_{t_s} = \sum_{m =1}^L \sum_{n =1}^L
\mu^{mn} \eta^{ip} F_n^{(s)} (u) {\pa^2 \psi_m \over \pa u^p \pa u^k}
u^k_x.
 \label{ierarsl}
\ee
All the flows in the hierarchy (\ref{ierar}) are commuting
integrable bi-Hamiltonian systems of hydrodynamic type
with an infinite family of local integrals in involution
with respect to both Poisson brackets:
\be
u^i_{t_s} =
M^{ij}_1 {\delta H_s
\over \delta u^j (x)} = \{ u^i (x), H_s \}_1, \ \ \ \
H_s = \int h_s (u(x)) dx,
\ee
\be
u^i_{t_s} = M^{ij}_2  {\delta H_{s+1}
\over \delta u^j (x)} = \{ u^i (x), H_{s+1} \}_2, \ \ \ \
H_{s+1} = \int h_{s+1} (u(x)) dx,
\ee
\be
 \{ H_p, H_r \}_1 = 0, \ \ \ \
\{ H_p, H_r \}_2 = 0,
\ee
and the densities $h_s (u(x))$ of the Hamiltonians
are related by the recursion relations (\ref{fsn}), (\ref{hsn}),
which are always solvable in our case.

\section{Locality and integrability of Hamiltonian
systems with \\nonlocal Poisson brackets of hydrodynamic type} \label{lok}

Let the functions
$\psi_n (u),$ $1 \leq n \leq L,$ be a solution
of the integrable system (\ref{g3}), (\ref{r3})
of nonlinear equations describing all flat torsionless
submanifolds in pseudo-Euclidean spaces; in particular, in this case
the nonlocal operator $M^{ij}_1$ given by the formula (\ref{nonlmd})
is Hamiltonian and compatible with the constant Hamiltonian operator
$M^{ij}_2$ {\rm (\ref{nonlme})}.

Consider the Hamiltonian system
\be
u^i_t =
M^{ij}_1 {\delta H
\over \delta u^j (x)} = \{ u^i (x), H \}_1   \label{arb}
\ee
with an arbitrary Hamiltonian of hydrodynamic type
\be
H = \int h (u(x)) dx. \label{arb2}
\ee
Ferapontov proved in \cite{1f} that
a Hamiltonian system with a nonlocal Hamiltonian
operator of hydrodynamic type
(\ref{nonle}) and with a Hamiltonian of hydrodynamic type (\ref{arb2})
is local if and only if the Hamiltonian
is an integral of all the structural flows of the nonlocal
Hamiltonian operator.
This statement is also true for Hamiltonian operators
of the form (\ref{nonlmb}), and moreover, it is always true
for any weakly nonlocal Hamiltonian
operators (see \cite{4}). We prove that
for the nonlocal Hamiltonian operators
$M^{ij}_1$ (\ref{nonlmd}) given by solutions of
the integrable system (\ref{g3}), (\ref{r3})
this condition on the Hamiltonians is sufficient
for integrability, i.e., all the corresponding local
Hamiltonian systems (\ref{arb}), (\ref{arb2}) are
integrable bi-Hamiltonian systems.

\bl  {\rm \cite{3m}}
The system {\rm (\ref{arb}), (\ref{arb2})} is local
if and only if the density $h (u(x))$ of the Hamiltonian
satisfies the linear equations
\be
{\pa^2 \psi_n \over
\pa u^j \pa u^s} \eta^{jr}  {\pa^2 h
\over \pa u^r \pa u^p}  =
 {\pa^2 \psi_n \over
\pa u^j \pa u^p} \eta^{jr}  {\pa^2 h
\over \pa u^r \pa u^s},  \ \ \  1 \leq n \leq L. \label{locint}
\ee
\el

Consider a system (\ref{arb}), (\ref{arb2}):
\be
u^i_t = M^{ij}_1 {\delta H
\over \delta u^j (x)} =
\sum_{m =1}^L \sum_{n =1}^L
\mu^{mn} \eta^{ip} {\pa^2 \psi_m \over \pa u^p \pa u^k}
u^k_x \left ( {d \over dx} \right )^{-1} \left ( \eta^{jr}
{\pa^2 \psi_n \over \pa u^r \pa u^s} u^s_x {\pa h
\over \pa u^j} \right ). \label{hamsyst}
\ee

The system (\ref{hamsyst}) is local if and only if
there locally exist functions $P_n (u),$ $1 \leq n \leq L,$
such that
\be
{\pa^2 \psi_n \over
\pa u^j \pa u^s} \eta^{jr}  {\pa h
\over \pa u^r} = {\pa P_n \over \pa u^s}, \label{pn}
\ee
i.e., if and only if the consistency relation
\be
{\pa \over \pa u^p} \left ( {\pa^2 \psi_n \over
\pa u^j \pa u^s} \eta^{jr}  {\pa h
\over \pa u^r} \right ) =
{\pa \over \pa u^s} \left ( {\pa^2 \psi_n \over
\pa u^j \pa u^p} \eta^{jr}  {\pa h
\over \pa u^r} \right ) \label{sovm}
\ee
is satisfied.
Then the system (\ref{hamsyst}) takes a local form
\be
u^i_t = M^{ij}_1 {\delta H
\over \delta u^j (x)} =
\sum_{m =1}^L \sum_{n =1}^L
\mu^{mn} \eta^{ip} {\pa^2 \psi_m \over \pa u^p \pa u^k} P_n (u)
u^k_x. \label{hamsyst2}
\ee

The consistency relation (\ref{sovm}) is equivalent to
the linear equations (\ref{locint}).

\bt  {\rm \cite{3m}}
If the functions
$\psi_n (u),$ $1 \leq n \leq L,$ are a solution
of the integrable system {\rm (\ref{g3}), (\ref{r3})} and
the corresponding Hamiltonian system {\rm (\ref{arb}), (\ref{arb2})}
is local, i.e., the density $h(u(x))$ of the Hamiltonian
satisfies the linear equations
{\rm (\ref{locint})}, then this Hamiltonian system is integrable
and bi-Hamiltonian.
\et

{\it Proof.}
In this case the system (\ref{arb}),
(\ref{arb2}) takes the form  (\ref{hamsyst2}),
(\ref{pn}).
Let us prove that there always locally exists a function
$f (u)$ such that
\be
\sum_{m =1}^L \sum_{n =1}^L
\mu^{mn} {\pa^2 \psi_m \over \pa u^j \pa u^k} P_n (u) =
{\pa^2 f \over \pa u^j \pa u^k}.
\ee

Indeed, we have
\bea
&&
{\pa \over \pa u^i} \left ( \sum_{m =1}^L \sum_{n =1}^L
\mu^{mn} {\pa^2 \psi_m \over \pa u^j \pa u^k} P_n (u) \right ) =
\sum_{m =1}^L \sum_{n =1}^L
\mu^{mn} {\pa^3 \psi_m \over \pa u^i \pa u^j \pa u^k} P_n (u) + \nn\\
&&
+ \sum_{m =1}^L \sum_{n =1}^L
\mu^{mn} {\pa^2 \psi_m \over \pa u^j \pa u^k} {\pa P_n  \over
\pa u^i} = \sum_{m =1}^L \sum_{n =1}^L
\mu^{mn} {\pa^3 \psi_m \over \pa u^i \pa u^j \pa u^k} P_n (u)
+ \nn\\
&&
+ \sum_{m =1}^L \sum_{n =1}^L
\mu^{mn} {\pa^2 \psi_m \over \pa u^j \pa u^k}
{\pa^2 \psi_n \over
\pa u^i \pa u^p} \eta^{pr}  {\pa h
\over \pa u^r} = \sum_{m =1}^L \sum_{n =1}^L
\mu^{mn} {\pa^3 \psi_m \over \pa u^i \pa u^j \pa u^k} P_n (u)
+ \nn\\
&&
+ \sum_{m =1}^L \sum_{n =1}^L
\mu^{mn} {\pa^2 \psi_m \over \pa u^j \pa u^i} {\pa^2 \psi_n  \over
\pa u^k \pa u^p} \eta^{pr}  {\pa h
\over \pa u^r},
\eea
where we have used the relation (\ref{g3}).
Consequently, by virtue of symmetry with
respect to the indices $i$ and $j$, we obtain
\be
{\pa \over \pa u^i} \left ( \sum_{m =1}^L \sum_{n =1}^L
\mu^{mn} {\pa^2 \psi_m \over \pa u^j \pa u^k} P_n (u) \right ) =
{\pa \over \pa u^j} \left ( \sum_{m =1}^L \sum_{n =1}^L
\mu^{mn} {\pa^2 \psi_m \over \pa u^i \pa u^k} P_n (u) \right ),
\ee
i.e., there locally exist functions
$b_k (u),$ $1 \leq k \leq N,$ such that
\be
 \sum_{m =1}^L \sum_{n =1}^L
\mu^{mn} {\pa^2 \psi_m \over \pa u^j \pa u^k} P_n (u) =
{\pa b_k \over \pa u^j}.
\ee
By virtue of symmetry with respect to the indices
$j$ and $k$, we obtain
\be
{\pa b_k \over \pa u^j} = {\pa b_j \over \pa u^k},
\ee
i.e., there locally exists a function $f (u)$
such that
\be
b_k (u) = {\pa f \over \pa u^k}.
\ee
Thus we have
\be
 \sum_{m =1}^L \sum_{n =1}^L
\mu^{mn} {\pa^2 \psi_m \over \pa u^j \pa u^k} P_n (u) =
{\pa b_k \over \pa u^j} = {\pa^2 f \over \pa u^j \pa u^k}.
\ee
Consequently, the system (\ref{arb}), (\ref{arb2})
in the case under consideration can be presented in the form
\be
u^i_t = \eta^{ij} {\pa^2 f \over \pa u^j \pa u^k} u^k_x =
M^{ij}_2 {\delta F
\over \delta u^j (x)} = \{ u^i (x), F \}_2,
 \ \ \ \ F = \int f (u) dx, \label{arb3}
\ee
i.e., it is an integrable bi-Hamiltonian system with
the compatible Hamiltonian operators $M^{ij}_1$
(\ref{nonlmd}) and $M^{ij}_2$ (\ref{nonlme}).

Thus we have an integrable description
(\ref{g3}), (\ref{r3}) and (\ref{locint}) of the class of local
Hamiltonian systems of the form (\ref{arb}), (\ref{arb2})
(we have proved that each of these local Hamiltonian systems is
integrable and bi-Hamiltonian).

\section{Systems of integrals in involution
generated by arbitrary
flat \\torsionless submanifolds
in pseudo-Euclidean spaces} \label{invol}

The nonlinear equations of the form
(\ref{r3}) and (\ref{locint})
are of independent interest. They play
an important role and have a
very natural interpretation.

\bl {\rm \cite{3m}}
The nonlinear equations {\rm (\ref{r3})} are equivalent
to the condition that the integrals
\be
\Psi_n = \int \psi_n (u(x)) dx, \ \ 1 \leq n \leq L,  \label{psin}
\ee
are in involution with respect to the Poisson bracket
defined by the constant Hamiltonian operator
$M^{ij}_2$ {\rm (\ref{nonlme})}, i.e., the condition
\be
\{ \Psi_n, \Psi_m \}_2 = 0, \ \ \ 1 \leq n, m \leq L.
\ee
\el

{\it Proof.}
Indeed, we have
\be
\{ \Psi_n, \Psi_m \}_2 = \int {\pa \psi_n \over \pa u^i}
\eta^{ij} {d \over d x} {\pa \psi_m \over \pa u^j} dx =
\int {\pa \psi_n \over \pa u^i}
\eta^{ij} {\pa^2 \psi_m \over \pa u^j \pa u^k} u^k_x dx.
\ee
Consequently, the integrals are in involution, i.e.,
\be
\{ \Psi_n, \Psi_m \}_2 = 0,
\ee
if and only if there exists a function
$S_{nm} (u)$ such that
\be
{\pa \psi_n \over \pa u^i}
\eta^{ij} {\pa^2 \psi_m \over \pa u^j \pa u^k} =
{\pa S_{nm} \over \pa u^k},
\ee
i.e., if and only if the consistency relation
\be
{\pa \over \pa u^l} \left ( {\pa \psi_n \over \pa u^i}
\eta^{ij} {\pa^2 \psi_m \over \pa u^j \pa u^k} \right ) =
{\pa \over \pa u^k} \left ( {\pa \psi_n \over \pa u^i}
\eta^{ij} {\pa^2 \psi_m \over \pa u^j \pa u^l} \right ) \label{sovm2}
\ee
is satisfied.
The consistency relation (\ref{sovm2}) is equivalent to
the equations (\ref{r3}).

Likewise,
the equations (\ref{locint}) are equivalent to the condition
\be
\{ \Psi_n, H \}_2 = 0, \ \ \ \ H = \int h (u(x)) dx.
\ee

We note that the equations (\ref{non1}) are
equivalent to the condition that $L$
integrals are in involution with respect to
an arbitrary Dubrovin--Novikov bracket (a nondegenerate
local Poisson bracket of hydrodynamic type).

\bt
If the functions
$\psi_n (u),$ $1 \leq n \leq L,$ are a solution
of the integrable system {\rm (\ref{g3}), (\ref{r3})}
of nonlinear equations describing all flat torsionless
submanifolds in pseudo-Euclidean spaces, then
the integrals $\Psi_n = \int \psi_n (u (x)) dx$ {\rm (\ref{psin})}
are in involution with respect to both the Poisson brackets
given by the nonlocal Hamiltonian operator $M^{ij}_1$
{\rm (\ref{nonlmd})} and the constant Hamiltonian operator
$M^{ij}_2$ {\rm (\ref{nonlme})}, and therefore the integrals are
in involution with respect to the corresponding pencil of
Poisson brackets{\rm :}
\be
\{ \Psi_n, \Psi_m \}_1 = 0, \ \ \
\{ \Psi_n, \Psi_m \}_2 = 0, \ \ \ 1 \leq n, m \leq L.
\ee
\et

\bt {\rm \cite{3m}}
If the functions
$\psi_n (u),$ $1 \leq n \leq L,$ are a solution
of the integrable system {\rm (\ref{g3}), (\ref{r3})}, then
the corresponding Hamiltonian system {\rm (\ref{arb}), (\ref{arb2}) }
is local if and only if it is generated by a family
of $L + 1$ integrals in involution with respect to
the Poisson bracket defined by the constant Hamiltonian operator
$M^{ij}_2$ {\rm (\ref{nonlme});} namely,
\be
\Psi_n = \int \psi_n (u(x)) dx,
1 \leq n \leq L, \ H = \int h (u(x)) dx, \ \
\{ \Psi_n, \Psi_m \}_2 = 0, \ \{ \Psi_n, H \}_2 = 0,
\ \ 1 \leq n, m \leq L.
\ee
Moreover, in this case
the system {\rm (\ref{arb}), (\ref{arb2})} is
an integrable bi-Hamiltonian system and
\be
\{ \Psi_n, \Psi_m \}_1 = 0, \ \{ \Psi_n, H \}_1 = 0,
\ \ 1 \leq n, m \leq L.
\ee
\et

\section{Systems of integrals in involution
generated by the \\associativity equations} \label{invol2}

An important special class of integrals in involution
is generated by the associativity equations
of two-dimensional topological quantum field theory
(the WDVV equations).

\bt {\rm \cite{3m}, \cite{5m}}
A function $\Phi (u^1, ..., u^N)$ generates a family
of $N$ integrals in involution with respect to
the Poisson bracket defined by the constant Hamiltonian operator
$M^{ij}_2$ {\rm (\ref{nonlme})}, namely, integrals whose densities are
the first-order partial derivatives of the function
{\rm (}the potential\/{\rm )} $\Phi (u)$
\be
I_n = \int {\pa \Phi \over \pa u^n} (u(x)) dx, \ \
\{ I_n, I_m \}_2 = 0, \ \  1 \leq n, m \leq N,  \label{intf}
\ee
if and only if the function $\Phi (u)$ is a solution
of the associativity equations {\rm (\ref{ass1a})} of
two-dimensional topological quantum field theory
{\rm (}the WDVV equations{\rm )}.
\et

Such special families of integrals in involution (whose densities are
the first-order partial derivatives of a potential $\Phi (u)$ and
the potential $\Phi (u)$ itself) with respect to
special nonlocal
Poisson brackets and with respect to pencils of
compatible Poisson brackets
generated by the associativity equations
will be also considered further in Section \ref{invol3}.
Note that we have constructed an infinite
family of integrals in involution with respect to both
local and nonlocal Poisson brackets generated by an arbitrary
flat torsionless submanifold in a pseudo-Euclidean space,
in particular, by an arbitrary Frobenius manifold, in Section \ref{hier5}.

\section{Associativity equations and nonlocal
Poisson brackets of \\hydrodynamic type} \label{sass}

Since the associativity equations (\ref{ass1a}) are
a natural reduction (see Section \ref{section3})
of the integrable system
{\rm (\ref{g3}), (\ref{r3})}
of nonlinear equations describing all flat torsionless
submanifolds in pseudo-Euclidean spaces and generating
all nonlocal Hamiltonian operators
of hydrodynamic type with flat metrics,
each solution $\Phi (u)$ of the associativity equations
(\ref{ass1a}), which are known to be
consistent and integrable by the inverse scattering method
and possess a rich set of nontrivial solutions
(see \cite{1d}), defines a nonlocal Hamiltonian operator
of hydrodynamic type with a flat metric
\be
L^{ij} = \eta^{ij} {d \over dx} +
\sum_{m=1}^N \sum_{n =1}^N
\eta^{mn} \eta^{ip} \eta^{jr}
{\pa^3 \Phi \over \pa u^p \pa u^m \pa u^k}
u^k_x \left ( {d \over dx} \right )^{-1} \circ
{\pa^3 \Phi \over \pa u^r \pa u^n \pa u^s} u^s_x \label{nonl7}
\ee
and even a pencil of compatible
Hamiltonian operators
\be
L^{ij}_{\lambda_1, \lambda_2} = \lambda_1
\eta^{ij} {d \over dx} +
\lambda_2 \sum_{m=1}^N \sum_{n =1}^N
\eta^{mn} \eta^{ip} \eta^{jr}
{\pa^3 \Phi \over \pa u^p \pa u^m \pa u^k}
u^k_x \left ( {d \over dx} \right )^{-1} \circ
{\pa^3 \Phi \over \pa u^r \pa u^n \pa u^s} u^s_x, \label{nonl7a}
\ee
where $\lambda_1$ and $\lambda_2$ are arbitrary constants.
In particular, if $\Phi (u)$ is an arbitrary
solution of the associativity equations
(\ref{ass1a}), then the operator
\be
L^{ij}_{0, 1} = \sum_{m=1}^N \sum_{n =1}^N
\eta^{mn} \eta^{ip} \eta^{jr}
{\pa^3 \Phi \over \pa u^p \pa u^m \pa u^k}
u^k_x \left ( {d \over dx} \right )^{-1} \circ
{\pa^3 \Phi \over \pa u^r \pa u^n \pa u^s} u^s_x \label{nonl7b}
\ee
is a Hamiltonian operator compatible with the
constant Hamiltonian operator
\be
L^{ij}_{1, 0} =
\eta^{ij} {d \over dx}.  \label{nonl7c}
\ee
The converse is also true.
\bt
The nonlocal operator $L^{ij}_{0, 1}$ {\rm (\ref{nonl7b})}
is Hamiltonian if and only if the function $\Phi (u)$ is a
solution of the associativity equations {\rm (\ref{ass1a})}.
\et
Therefore, for each solution of the associativity equations
(\ref{ass1a}) (in particular, for each Frobenius manifold)
we obtain the corresponding natural pencil of compatible
Poisson structures (local and nonlocal) and the
corresponding natural integrable hierarchies
(see Section \ref{hier5}).

Thus, for each Frobenius manifold
there are a very natural nonlocal
Hamiltonian operator of the form (\ref{nonl7}),
a pencil of compatible Hamiltonian operators
(\ref{nonl7a}) and very natural
integrable hierarchies connected to the Frobenius manifold.

We have considered the nonlocal Hamiltonian operators of
the form (\ref{nonlm}) with flat metrics and came to
the associativity equations defining the affinors of such operators.
A statement that is in some sense the converse is also true,
namely, if all the affinors $w_n (u)$
of a nonlocal Hamiltonian operator (\ref{nonlm})
with $L = N$ are defined by an arbitrary solution
$\Phi (u)$ of the associativity equations (\ref{ass1a})
by the formula
$$
(w_n)^i_j (u) = \zeta^{is} \xi^r_j {\pa^3 \Phi \over
\pa u^n \pa u^s \pa u^r},
$$
where $\zeta^{is},$ $\xi^r_j$ are arbitrary nondegenerate
constant matrices, then the metric of this Hamiltonian
operator must be flat. But, in general, it is not necessarily that
this metric will be constant in the local coordinates under consideration.

The structural flows (see \cite{1f}, \cite{4})
of the nonlocal Hamiltonian operator (\ref{nonl7})
have the form:
\be
u^i_{t_n} = \eta^{is} {\pa^3 \Phi \over \pa u^s \pa u^n \pa u^k}
u^k_x.  \label{nonl8}
\ee
These systems are integrable bi-Hamiltonian systems of
hydrodynamic type and coincide with the primary part of
the Dubrovin hierarchy constructed by any solution
of the associativity equations in \cite{1d}.
The condition of commutation for
the structural flows (\ref{nonl8})
is also equivalent to the associativity equations (\ref{ass1a}).

\bt
For an arbitrary solution
$\Phi (u)$ of the associativity equations {\rm (\ref{ass1a})},
each structural flow {\rm (\ref{nonl8})} generates an integrable
hierarchy of hydrodynamic type
with the recursion
operator given by the compatible Hamiltonian operators
{\rm (\ref{nonl7b})} and {\rm (\ref{nonl7c})}
{\rm (}see the recursion operator {\rm (\ref{rek})} in a more general
case{\rm )}{\rm ;} each of these integrable
hierarchies is local and bi-Hamiltonian
with respect to the compatible Hamiltonian operators
{\rm (\ref{nonl7b})} and {\rm (\ref{nonl7c})}.
\et

A great number of concrete examples of Frobenius
manifolds and solutions of the associativity equations
is given in Dubrovin's paper \cite{1d}.
Consider here only one simple example from
\cite{1d} as an illustration.
Let $N = 3$ and the metric $\eta_{ij}$ be antidiagonal
\be
(\eta_{ij}) =
\left ( \begin{array} {ccc} 0&0&1\\
0&1&0\\
1&0&0
\end{array} \right ),
\ee
and the function $\Phi (u)$ has the form
$$
\Phi (u) = {1 \over 2} (u^1)^2 u^3 +
{1 \over 2} u^1 (u^2)^2 + f (u^2, u^3).
$$
In this case $e_1$ is the unit in the Frobenius algebra
(\ref{frob}), (\ref{form}), and the associativity equations
(\ref{ass1a}) for the function $\Phi (u)$ are equivalent to
the following remarkable integrable Dubrovin equation for the function
$f (u^2, u^3)$:
\be
{\pa^3 f \over \pa (u^3)^3} = \left (
{\pa^3 f \over \pa (u^2)^2 \pa u^3} \right )^2 -
{\pa^3 f \over \pa (u^2)^3} {\pa^3 f \over \pa u^2 \pa (u^3)^2}.
\label{f}
\ee
This equation is connected to quantum cohomology of
projective plane and classical problems of
enumerative geometry (see \cite{6}). In particular,
all nontrivial polynomial solutions of the equation (\ref{f})
that satisfy the requirement of the quasihomogeneity
and locally define a structure of
Frobenius manifold are described by Dubrovin in \cite{1d}:
\be
f = {1 \over 4} (u^2)^2 (u^3)^2 + {1 \over 60} (u^3)^5,\ \ \
f = {1 \over 6} (u^2)^3 u^3 + {1 \over 6} (u^2)^2 (u^3)^3 +
{1 \over 210} (u^3)^7, \label{sol1}
\ee
\be
f = {1 \over 6} (u^2)^3 (u^3)^2 + {1 \over 20} (u^2)^2 (u^3)^5 +
{1 \over 3960} (u^3)^{11}. \label{sol2}
\ee

As is shown by the author in \cite{7} (see also \cite{8}),
the equation (\ref{f}) is equivalent to the integrable
nondiagonalizable homogeneous system of hydrodynamic type
\be
\left ( \begin{array} {c} a\\ b\\ c
\end{array} \right )_{u^3} =
\left ( \begin{array} {ccc} 0 & 1 & 0\\  0 & 0 & 1\\
- c & 2b & - a
\end{array} \right )  \left ( \begin{array} {c} a\\ b\\ c
\end{array} \right )_{u^2},
\ee
\be
a = {\pa^3 f \over \pa (u^2)^3}, \ \ \
b = {\pa^3 f \over \pa (u^2)^2 \pa u^3},\ \ \
c = {\pa^3 f \over \pa u^2 \pa (u^3)^2}.
\ee
In this case the affinors of the nonlocal Hamiltonian
operator
(\ref{nonl7}) have the form:
\be
(w_1)^i_j (u) = \delta^i_j, \ \ \
(w_2)^i_j (u) =
\left ( \begin{array} {ccc} 0 & b & c\\  1 & a & b\\
0 & 1 & 0
\end{array} \right ),\ \ \
(w_3)^i_j (u) = \left ( \begin{array} {ccc}
0 & c & b^2 - a c\\  0 & b & c\\
1 & 0 & 0
\end{array} \right ).
\ee

For concrete solutions of the associativity equation (\ref{f}),
in particular, for (\ref{sol1}) and (\ref{sol2}),
the corresponding linear systems (\ref{rf1}),
(\ref{rf2}) giving explicit realizations of
the corresponding Frobenius manifolds as potential
flat torsionless submanifolds in pseudo-Euclidean spaces
can be solved in special functions;
we shall give the explicit realizations
in a separate paper.

\section{Associativity equations and special integrals in involution
with respect to nonlocal Poisson
brackets of hydrodynamic type and Poisson pencils} \label{invol3}

If the function $\Phi (u^1, \ldots, u^N)$ is an arbitrary
solution of the associativity equations
(\ref{ass1a}), then the operator
$L^{ij}_{0, 1}$ (\ref{nonl7b})
is a Hamiltonian operator, and we can consider the corresponding
Poisson bracket and integrals in involution with respect to
this Poisson bracket.

\bt {\rm \cite{5m}}
If the function $\Phi (u^1, \ldots, u^N)$ satisfies
the associativity equations
{\rm (\ref{ass1a})}, then the functionals
$I_n = \int (\partial \Phi/ \partial u^n) dx$, $n = 1, \ldots, N,$
{\rm (\ref{intf})} are in involution with respect to the
Poisson bracket given by the nonlocal Hamiltonian operator
$L^{ij}_{0, 1}$ {\rm (\ref{nonl7b})}, i.e.,
\be
\ \ \{ I_n, I_m \}_1 = 0,
\ \ \ \ n, m = 1, \ldots, N, \label{7}
\ee
where
\be
\ \ \{u^i (x), u^j (y) \}_1 =
\eta^{mn} \eta^{ip} \eta^{jr}
{\partial^3 \Phi \over \partial u^p \partial u^m \partial u^k}
u^k_x \left ( {d \over dx} \right )^{- 1}
\left ( {\partial^3 \Phi \over \partial u^r \partial u^n \partial u^s}
u^s_x  \delta (x - y) \right ). \label{8}
\ee
\et
The Poisson brackets
$\{u^i (x), u^j (y)\}_2 = \eta^{ij} \delta^{\prime} (x - y)$ given
by the constant Hamiltonian operator
$L^{ij}_{1, 0}$ (\ref{nonl7c}) and $\{ u^i (x), u^j (y) \}_1$
given by the nonlocal Hamiltonian operator
$L^{ij}_{0, 1}$ {\rm (\ref{nonl7b})}
are compatible and form a pencil
$\lambda_1 \{ u^i (x), u^j (y) \}_1 +
\lambda_2 \{ u^i (x), u^j (y) \}_2$ of Poisson brackets,
where $\lambda_1$ and $\lambda_2$ are arbitrary constants, so that
for any constants $\lambda_1$ and $\lambda_2$ the bracket
$$\{ u^i (x), u^j (y) \}_{\lambda_1, \lambda_2}
= \lambda_1 \{ u^i (x), u^j (y) \}_1 +
\lambda_2 \{ u^i (x), u^j (y) \}_2$$
is a Poisson bracket.

\bc {\rm \cite{5m}}
If the function $\Phi (u^1, \ldots, u^N)$ satisfies
the associativity equations
{\rm (\ref{ass1a})}, then the functionals
$I_n = \int (\partial \Phi/ \partial u^n) dx$, $n = 1, \ldots, N,$
{\rm (\ref{intf})} are in involution with respect to the pencil of
Poisson brackets
$\{ u^i (x), u^j (y) \}_{\lambda_1, \lambda_2}
= \lambda_1 \{ u^i (x), u^j (y) \}_1 +
\lambda_2 \{ u^i (x), u^j (y) \}_2$,
where $\lambda_1$ and $\lambda_2$  are arbitrary constants.
\ec

We consider also the functional
$F$ whose density is the potential
$\Phi (u^1 (x), \ldots, u^N (x))$ itself:
\be
 F =
\int  \Phi (u^1 (x), \ldots, u^N (x)) dx. \label{9}
\ee

\bt {\rm \cite{5m}}
A function $\Phi (u^1, ..., u^N)$ generates a family
of $N + 1$ integrals in involution with respect to
the constant Poisson bracket
$\{u^i (x), u^j (y)\}_2 = \eta^{ij} \delta^{\prime} (x - y)$,
namely, the functional $F$ {\rm (\ref{9})} and
the functionals
$I_n = \int (\partial \Phi/ \partial u^n) dx$, $n = 1, \ldots, N,$
{\rm (\ref{intf})}, $\{ I_n, I_m \}_2 = 0,$ $1 \leq n, m \leq N,$
$\{ I_n, F \}_2 = 0,$ $1 \leq n \leq N,$
if and only if the function
$\Phi (u^1, \ldots, u^N)$ satisfies the equations
\be
{\partial^2 \Phi \over \partial u^k \partial u^i}
\eta^{ij}
{\partial^3 \Phi
\over \partial u^j \partial u^n \partial u^l} =
{\partial^2 \Phi \over \partial u^l \partial u^i}
\eta^{ij}
{\partial^3 \Phi \over \partial u^j \partial u^n \partial u^k}. \label{10}
\ee
\et

The equations (\ref{10}) have arised in a different context
in the author's papers \cite{6m}--\cite{8m}
and play an important role in the theory of compatible Poisson
brackets of hydrodynamic type, the theory of the associativity equations
and the theory of Frobenius manifolds.

\bt {\rm \cite{5m}}
If the function $\Phi (u^1, \ldots, u^N)$ satisfies
the equations
{\rm (\ref{10})}, then the functional $F$ {\rm (\ref{9})}
and the functionals
$I_n = \int (\partial \Phi/ \partial u^n) dx$, $n = 1, \ldots, N,$
{\rm (\ref{intf})} are in involution with respect to the
Poisson bracket $\{ u^i (x), u^j (y) \}_1$
given by the nonlocal Hamiltonian operator
$L^{ij}_{0, 1}$ {\rm (\ref{nonl7b})}.
\et

\bc {\rm \cite{5m}}
If the function $\Phi (u^1, \ldots, u^N)$ satisfies
the equations
{\rm (\ref{10})}, then the functional $F$ {\rm (\ref{9})}
and the functionals
$I_n = \int (\partial \Phi/ \partial u^n) dx$, $n = 1, \ldots, N,$
{\rm (\ref{intf})} are in involution with respect to the
pencil of
Poisson brackets
$\{ u^i (x), u^j (y) \}_{\lambda_1, \lambda_2}
= \lambda_1 \{ u^i (x), u^j (y) \}_1 +
\lambda_2 \{ u^i (x), u^j (y) \}_2$,
where $\lambda_1$ and $\lambda_2$  are arbitrary constants.
\ec



\end{document}